\documentclass[12pt]{article}
\usepackage[english]{babel}

\usepackage{amsmath, amsthm}
\usepackage{amscd}
\usepackage{amsfonts}
\usepackage{amssymb}
\usepackage{xypic}
\usepackage{tikz-cd}
\usepackage{mathabx}
\usepackage{mathtools}
\usepackage{graphicx}
\usepackage{hyperref}

\usepackage[left=2cm,right=2cm,
    top=2cm,bottom=2cm,bindingoffset=0cm]{geometry}

\title{THE WEIL CONJECTURES I}
\author{By Pierre Deligne (translated by Evgeny Goncharov) \\ \href{mailto:eg555@cam.ac.uk}{eg555@cam.ac.uk}\footnote{Cambridge University, Center for Mathematical Sciences, Wilberforce Road, Cambridge, UK}, \href{mailto:eagoncharov@edu.hse.ru}{eagoncharov@edu.hse.ru}\footnote{National Research University Higher School of Economics (NRU HSE), Usacheva 6, Moscow, Russia.}}
\date{}

\newtheorem*{fact1}{Proposition \textbf{(5.3)}}
\newtheorem*{fact2}{Proposition \textbf{(6.6)}}
\newtheorem*{fact3}{Proposition \textbf{(6.8)}}

\newtheorem*{lemma}{Lemma \textbf{(1.7)}}
\newtheorem*{lemma1}{Lemma \textbf{(3.3)}}
\newtheorem*{lemma2}{Lemma \textbf{(3.4)}}
\newtheorem*{lemma3}{Lemma \textbf{(3.5)}}
\newtheorem*{lemma4}{Lemma \textbf{(3.6)}}
\newtheorem*{lemma5}{Lemma \textbf{(5.11)}}
\newtheorem*{lemma6}{Lemma \textbf{(6.4)}}
\newtheorem*{lemma7}{Lemma \textbf{(6.7)}}
\newtheorem*{lemma8}{Lemma \textbf{(6.11)}}
\newtheorem*{lemma9}{Lemma \textbf{(6.12)}}
\newtheorem*{lemma10}{Lemma \textbf{(7.1)}}
\newtheorem*{lemma11}{Lemma \textbf{(7.2)}}
\newtheorem*{lemma12}{Lemma \textbf{(8.5)}}
\newtheorem*{lemma13}{Lemma \textbf{(8.7)}}
\newtheorem*{lemma14}{Lemma \textbf{(8.8)}}
\newtheorem*{lemma15}{Lemma \textbf{(8.12)}}
\newtheorem*{lemma16}{Lemma \textbf{(8.13)}}
\newtheorem*{thrm1}{Theorem \textbf{(1.6)}}
\newtheorem*{thrm2}{Theorem \textbf{(2.3)} (Poincare duality)}
\newtheorem*{thrm3}{Theorem \textbf{(2.8)}}
\newtheorem*{thrm4}{Theorem \textbf{(2.12)}}
\newtheorem*{thrm5}{Theorem \textbf{(3.2)}}
\newtheorem*{thrm6}{Theorem \textbf{(5.4)}}
\newtheorem*{thrm7}{Theorem \textbf{(5.10)} (Kajdan-Margulis)}
\newtheorem*{thrm8}{Theorem \textbf{(6.2)}}
\newtheorem*{thrm9}{Theorem \textbf{(8.1)}}
\newtheorem*{thrm10}{Theorem \textbf{(8.2)}}
\newtheorem*{thrm11}{Theorem \textbf{(8.4)}}
\newtheorem*{cor}{Corollary \textbf{(3.8)}}
\newtheorem*{cor1}{Corollary \textbf{(3.9)}}
\newtheorem*{cor2}{Corollary \textbf{(5.5)}}
\newtheorem*{cor3}{Corollary \textbf{(6.3)}}

\newtheorem*{stat}{Statement \textbf{(2.10)}}

\begin{document}

\maketitle

\textit{I attempted to write the full translation of this article to make the remarkable proof of Pierre Deligne available to a greater number of people. Overviews of the proofs can be found elsewhere. I especially recommend the notes of James Milne on Etale Cohomology that also contain a justification for the theory underlying this article and proofs of the results used by Deligne. The footnotes are mostly claims that some details appear in Milne, clarifications of some of the terminology or my personal struggles. I have also made a thorough overview of the proof together with more detailed explanations - \url{https://arxiv.org/abs/1807.10812}. Enjoy!}

\begin{abstract}
In this article I prove the Weil conjecture about the eigenvalues of Frobenius endomorphisms. The precise formulation is given in (1.6). I tried to make the demonstration as geometric and elementary as possible and included reminders: only the results of paragraphs 3, 6, 7 and 8 are original. 

In the article following this one I will give various refinements of the intermediate results and the applications, including the hard Lefschetz theorem (on the iterated cup products by the cohomology class of a hyperplane section).

The text faithfully follows from the six lectures given at Cambridge in July 1973. I thank N.Katz for allowing me to use his notes.
\end{abstract}

\tableofcontents

\section{The theory of Grothendieck: a cohomological interpretation of L-functions}
\textbf{(1.1)} Let $X$ be a scheme of finite type over $\mathbb{Z}$, $|X|$ be the set of closed points of $X$ and for $x \in |X|$ we denote by $N(x)$ the number of elements in the residue field $k(x)$ of $X$ at $x$. The Hasse-Weil zeta function of $X$ is $$\zeta_X(s)=\prod_{x \in |X|}(1-N(x)^{-s})^{-1} \ \ \textbf{(1.1.1)}$$ (this product converges absolutely for $Re(s)$ large enough). For $X=Spec(\mathbb{Z})$, $\zeta_X(s)$ is the Riemann zeta function. 

We will consider exclusively the case when $X$ is a scheme over a finite field $\mathbb{F}_q$. 

For $x \in |X|$ we will write $g_x$ instead of $N(x)$. Denoting $\deg(x)=[k(x) : \mathbb{F}_q]$ we have $q_x=q^{\deg(x)}$. It makes sense to introduce a new variable $t=q^{-s}$. Let $$Z(X; t)=\prod_{x \in |X|}(1-t^{\deg(x)})^{-1} \ \ \textbf{(1.1.2)};$$ this product converges for $|t|$ small enough and we have $$\zeta_X(s)=Z(X;q^{-s}) \ \ \textbf{(1.1.3)}$$

\textbf{(1.2)} Dwork (On the rationality of the zeta function of an algebraic variety, Amer. J. Math., 82, 1960, p. 631-648) and Grothendieck ([1] and SGA5) have demonstrated that $Z(X; t)$ is a rational function of $t$. 

For Grothendieck, this is a corollary of general results in $l$-adic cohomology (where $l$ is a prime number not equal to the characteristic $p$ of $\mathbb{F}_q$). These provide a cohomological interpretation of the zeros and poles of $Z(X; t)$, and a functional equation when $X$ is proper and smooth. The methods of Dwork are $p$-adic. For X a non-singular hypersurface in a projective space they also provided him with a cohomological interpretation of zeros and poles, and the functional equation. They inspired the crystalline theory of Grothendieck and Berthelot, which for X proper and smooth provides a $p$-adic cohomological interpretation of zeros and poles, and the functional equation. Based on Washnitzer ideas, Lubkin created a variant of this theory, valid only for X proper, smooth and liftable to characteristic $0$ (A $p$-adic proof of Weil's conjectures, Ann of Math, 87, 1968, pp. 125-255).

We will make essential use of Grothendieck's results and recall them below. 

\textbf{(1.3)} Let $X$ be an algebraic variety over an algebraically closed field $k$ of characteristic $p$, i.e. a separated scheme of finite type over $k$. We do not exclude the case $p=0$. For any prime number $l \neq p$, Grothendieck defined $l$-adic cohomology groups $H^i(X, \mathbb{Q}_l)$. He also defined cohomology groups with compact support $H_c^i(X, \mathbb{Q}_l)$. For $X$ proper the two coincide. $H_c^i(X, \mathbb{Q}_l)$ are vector spaces of finite dimension over $\mathbb{Q}_l$, zero for $i>2dim(X)$.

\textbf{(1.4)} Let $X_0$ be an algebraic variety over $\mathbb{F}_q$, $\mathbb{\bar F}_q$ the algebraic closure of $\mathbb{F}_q$ and $X$ the algebraic variety over $\mathbb{\bar F}_q$ obtained from $X_0$ by extension of scalars of $\mathbb{F}_q$ to $\mathbb{\bar F}_q$. In the language of Weil and Shimura we would express this situation by: "Let $X$ be an algebraic variety defined over $\mathbb{F}_q$". Let $F: X \to X$ be the Frobenius morphism; it sends a point with coordinates $x$ to the point with coordinates $x^q$; in other words, for $U_0$ a Zariski open subset of $X_0$, defining an open subset $U$ of $X$, we have $F^{-1}(U)=U$; for $x \in H^0(U_0, \mathcal{O})$ we have $F^*x=x^q$. Let us identify the set $|X|$ of closed points of $X$ with $X_0(\mathbb{\bar F}_q)$ (all the points $Hom_{\mathbb{F}_q}(Spec(\mathbb{\bar F}_q), X_0)$ of $X_0$ with coefficients in $\mathbb{\bar F}_q$) and let $\varphi \in Gal(\mathbb{\bar F}_q / \mathbb{F}_q)$ be the substitution of Frobenius: $\varphi(x)=x^q$. The action of $F$ on $|X|$ identifies with the action of $\varphi$ on $X_0(\mathbb{\bar F}_q)$. Then:

a) The set $X^F$ of closed points of $X$ fixed under $F$ is identified with the set $X_0(\mathbb{F}_q) \subset X_0(\mathbb{\bar F}_q)$ of points of $X$ defined over $\mathbb{F}_q$. This simply expresses the fact that for $x \in \mathbb{\bar F}_q$ we have $x \in \mathbb{F}_q \Leftrightarrow x^q=x$.

b) Similarly, the set $X^{F^n}$ of closed points of $X$ fixed under the $n$-th iteration of $F$ is identified with $X_0(\mathbb{F}_{q^n})$.

c) The set $|X|$ of closed points of $X$ is identified with the set $|X|_F$ of orbits of $F$ (or $\varphi$) on $|X|$. The degree $\deg(x)$ of $x \in |X_0|$ is the number of elements in the corresponding orbit.

d) From b) and c) we see that $$\#X^{F^n}=\#X_0(\mathbb{F}_{q^n})=\sum_{\deg(x)|n} \deg(x) \ \ \textbf{(1.4.1)}$$ (for $x \in |X_0|$ and $\deg(x)|n$, $x$ defines $\deg(x)$ points with coordinates in $\mathbb{F}_{q^n}$ all conjugate over $\mathbb{F}_q$).

\textbf{(1.5)} The morphism $F$ is finite, in particular, proper. Therefore, it induces morphisms $$F^*: H_c^i(X, \mathbb{Q}_l) \to H_c^i(X, \mathbb{Q}_l).$$ Grothendieck proved the formula of Lefschetz $$\#X^F=\sum_i (-1)^iTr(F^*, H_c^i(X, \mathbb{Q}_l));$$ the right side, that is a priori an $l$-adic number is an integer and is equal to the left side. We should note that such a formulation is only reasonable because $dF=0$, even at infinity ($X$ is not assumed to be proper); the relation $dF=0$ implies that fixed points of $F$ have multiplicity one.

We have a similar formula for the iterations of $F$: $$\#X^{F^n}=X_0(\mathbb{F}_{q^n})=\sum_i (-1)^i Tr(F^{*n}, H_c^i(X, \mathbb{Q}_l)) \ \ \textbf{(1.5.1)}$$

We take the logarithmic derivative of (1.1.2): $$t \frac{d}{dt}\log Z(X_0, t)=\frac{t\frac{d}{dt}Z(X_0, t)}{Z(X_0, t)}=\sum_{x \in |X_0|} -\frac{-\deg(x)t^{\deg(x)}}{1-t^{\deg(x)}}=$$ $$=\sum_{x \in |X_0|} \sum_{n >0} \deg(x)t^{n\deg(x)} \overset{(1.4.1)} = \sum_n X_0(\mathbb{F}_{q^n})t^n \ \ \textbf{(1.5.2)}$$

For $F$ an endomorphism of a vector space $V$ we have a formal series identity $$t \frac{d}{dt}\log(\det(1-Ft, V)^{-1})=\sum_{n > 0} Tr(F^n, V)t^n \ \ \textbf{(1.5.3)}$$ (check for $\dim V=1$ and observe that both sides are additive in V when we take short exact sequences). By substituting (1.5.1) into (1.5.2) and applying (1.5.3) one finds $$t \frac{d}{dt}\log Z(X_0, t)=\sum_i (-1)^i t \frac{d}{dt}\log \det (1-F^*t, H_c^i(X, \mathbb{Q}_l))^{-1},$$ or $$Z(X, t)=\prod_i \det (1-F^*t, H_c^i(X, \mathbb{Q}_l))^{(-1)^{i+1}} \ \ \textbf{(1.5.4)}$$ The right side is in $\mathbb{Q}_l(t)$. The formula implies that its Taylor expansion at $t=0$, a priori a formal series in $\mathbb{Q}_l [[t]]$ with constant coefficient one, is in $\mathbb{Z}[[t]]$ and is equal to the left side, also considered as a formal series in $t$. This formula is the Grothendiek's cohomological interpretation of the $Z$-function.

Our main result is the following:
\begin{thrm1}
Let $X_0$ be a projective nonsingular (= smooth) variety over $\mathbb{F}_q$. For each $i$, the characteristic polynomial $\det(1-F^*t, H^i(X, \mathbb{Q}_l))$ has integer coefficients independent of $l$ $(l \neq p)$. The complex roots $\alpha$ of this polynomial (complex conjugates of the eigenvalues of $F^*$) are of absolute value $|\alpha|=q^{\frac{i}{2}}$. 
\end{thrm1}

We show that (1.6) is a consequence of the following apparently weaker statement:

\begin{lemma}
For each $i$ and each $l \neq p$ the eigenvalues of the Frobenius endomorphism $F^*$ on $H^i(X, \mathbb{Q}_l)$ are algebraic numbers of absolute value $|\alpha|=q^{\frac{i}{2}}$.
\end{lemma}

\textit{Proof of (1.7) $\Rightarrow $(1.6)}: Let's look at $Z(X_0, t)$ as a formal series with constant term $1$, an element of $\mathbb{Z}[[t]]: Z(X_0, t)=\sum_n a_n t^n$. From (1.5.3), the image of $Z(X_0, t)$ in $\mathbb{Q}_l[[t]]$ is a Taylor expansion of a rational function. This means that for $N$ and $M$ large enough ($\geq$ the degrees of numerator and denominator) the Hankel determinants $$H_k=\det((a_{i+j+k})_{0 \leq i, j \leq M}) \ (k > N)$$ are zero. They vanish in $\mathbb{Q}_l$ if and only if they vanish in $\mathbb{Q}$; $Z(X_0, t)$ is a Taylor expansion of an element in $\mathbb{Q}(t)$. In other words, $$Z(X_0, t) \in \mathbb{Z}[[t]] \cap \mathbb{Q}_l(t) \subset \mathbb{Q}(t).$$

Let $Z(X_0, t)=\frac{P}{Q}$, with $P, Q \in \mathbb{Z}[t]$ coprime and with positive constant terms. According to a lemma of Fatou, since $Z(X_0, t)$ lies in $\mathbb{Z}[[t]]$ and has constant term $1$, the constant terms of $P$ and $Q$ are $1$\footnote{See the proof in James Milne's lectures on Etale Cohomology (Milne). Here and below - footnotes of the translator.}. Let $$P_i (t)=\det (1-F^*t, H^i (X, \mathbb{Q}_l)).$$ (1.7) implies that $P_i$ are coprime. The right hand side of (1.5.4) is therefore an irreducible fraction and $$ P(t)=\prod_{i \ odd} P_i(t)$$ $$ Q(t)=\prod_{i \ even} P_i(t).$$ Let $K$ be the subfield of the algebraic closure $\mathbb{\bar Q}_l$ of $\mathbb{Q}_l$ generated over $\mathbb{Q}$ by the roots of $R(t)=P(t)Q(t)$. The roots of $P_i(t)$ are the roots of $R(t)$ such that all their complex conjugates have absolute value $q^{-\frac{i}{2}}$. This set is stable under $Gal(K / \mathbb{Q})$. Therefore, $P_i(t)$ has rational coefficients. According to a lemma of Gauss (or because roots of $P_i$, being roots of $R$, are inverses of algebraic integers), it even has integer coefficients. The above description of the roots of $P_i(t)$ is independent of $l$, therefore, the polynomial $P_i(t)$ is also independent of $l$. 

The rest of the article is dedicated to the demonstration of (1.7). 

\textbf{(1.8)} The theory of Grothendieck provides cohomological interpretation not only of zeta functions but also of $L$-functions. The results are as follows.

\textbf{(1.9)} Let $X$ be an algebraic variety over a field $k$. For the definition of a constructible $\mathbb{Q}_l$-sheaf on $X$ consult SGA 5 VI\footnote{Or Milne.}. It suffices to say that: 

a) If $\mathcal{F}$ is a constructible $\mathbb{Q}_l$-sheaf, there exists a finite partition of $X$ into locally closed parts such that $\mathcal{F}|X_i$ is locally constant.

b) Assume that $X$ is connected and let $\bar x$ be a geometric point of $X$. For $\mathcal{F}$ locally constant, $\pi_1(X, \bar x)$ acts on the stalks $\mathcal{F}_{\bar x}$; the map $\mathcal{F} \to \mathcal{F}_{\bar x}$ defines an equivalence of categories (locally constant $\mathbb{Q}_l$-sheaves on $X$) $\to$ (continuous representations of $\pi_1(X, \bar x)$ on $\mathbb{Q}_l$ vector spaces of finite dimension). Such a representation in general does not factor through a finite quotient.

c) If $k=\mathbb{C}$, the constructible $\mathbb{Q}_l$-sheaves over $X$ are identified with the sheaves of $\mathbb{Q}_l$ vector spaces $\mathcal{F}$ on $X^{an}$ and their exists a finite partition of $X$ into Zariski locally closed parts and for each $i$ a local system\footnote{We will also call it a locally constant sheaf (french. constant tordu). This is an abuse of terminology. A sheaf $\mathcal{M}=(\mathcal{M}_n)$ of $\mathbb{Z}_l$-modules is called locally constant if each $\mathcal{M}_n$ is locally constant. It is not, in general, locally constant in the classical sense. Similar remarks apply elsewhere to $\mathbb{Z}_l$ and $\mathbb{Q}_l$-sheaves.} of free of finite type $\mathbb{Z}_l$-modules $\mathcal{F}_i$ on $X_i$ such that $$\mathcal{F}|X_i=\mathcal{F}_i \otimes_{\mathbb{Z}_l} \mathbb{Q}_l.$$

We will only consider constructible $\mathbb{Q}_l$-sheaves and call them just $\mathbb{Q}_l$-sheaves.

\textbf{(1.10)} Assume that $k$ is algebraically closed and let $\mathcal{F}$ be a $\mathbb{Q}_l$-sheaf on $X$. Grothendieck defined the $l$-adic cohomology groups $H^i(X, \mathcal{F})$ and $H_c^i(X, \mathcal{F})$. $H_c^i(X, \mathcal{F})$ are vector spaces of finite dimension over $\mathbb{Q}_l$, zero for $i > 2 \dim(X)$. For $k=\mathbb{C}, \  H^i(X, \mathcal{F})$ and $H_c^i(X, \mathcal{F})$ are the usual cohomology groups (resp. groups with compact support) of $X^{an}$ with coefficients in $\mathcal{F}$. 

\textbf{(1.11)} Let $X_0$ be an algebraic variety over $\mathbb{F}_q$, $X$ the corresponding variety over $\mathbb{\bar F}_q$ and $\mathcal{F}_0$ a sheaf of sets on $X_0$ (for the etale topology). We denote by $\mathcal{F}$ its inverse image on $X$. In addition to the Frobenius isomorphism $F: X \to X$, we have a canonical isomorphism $F^*: F^* \mathcal{F} \overset{\sim} \to \mathcal{F}$. Here is a description. We regard $\mathcal{F}_0$ as an etale space over $X_0$, i.e. we identify $\mathcal{F}_0$ with an algebraic space $[\mathcal{F}_0]$, equipped with an etale morphism $f: [\mathcal{F}_0] \to X_0$ such that $\mathcal{F}_0$ is the sheaf of local sections of $[\mathcal{F}_0]$. The similar etale space $[\mathcal{F}]$ over $X$ is obtained from $[\mathcal{F}_0]$ by extension of scalars. So we have a commutative diagram 
\[
\begin{tikzcd}
{[\mathcal{F}]} \arrow{r}{F} \arrow[swap]{d}{f} & {[\mathcal{F}]} \arrow{d}{f} \\
X \arrow{r}{F} & X
\end{tikzcd}
\]
and a morphism $[\mathcal{F}] \to X \times_{(F, X, f)} [\mathcal{F}]=[F^* \mathcal{F}]$, that is an isomorphism because $f$ is etale. The inverse of this isomorphism defines the isomorphism $F^* \mathcal{F} \overset{\sim} \to \mathcal{F}$ that we seek. 

This construction is generalized to $\mathbb{Q}_l$-sheaves.

\textbf{(1.12)} Let $X_0$ be an algebraic variety over $\mathbb{F}_q$, $\mathcal{F}_0$ a $\mathbb{Q}_l$-sheaf on $X_0$, $(X, \mathcal{F}$) is obtained by extension of scalars of $\mathbb{F}_q$ to $\mathbb{\bar F}_q$, $F: X \to X$ and $F^*: F^*\mathcal{F} \to \mathcal{F}$. Finite morphisms $F$ and $F^*$ define an endomorphism $$F^*: H_c^i(X, \mathcal{F}) \to H_c^i(X, F^*\mathcal{F}) \to H_c^i(X, \mathcal{F}).$$ For $x \in |X|$, $F^*$ defines a morphism $F_x^*:\mathcal{F}_{F(x)} \to \mathcal{F}_x$. For $x \in X^F$ it is an endomorphism of $\mathcal{F}_x$. Grothendieck proved the formula of Lefschetz $$\sum_{x \in X^F} Tr(F_x^*, \mathcal{F}_x)=\sum_i (-1)^i Tr(F^*, H_c^i(X, \mathcal{F})).$$ A similar formula holds for the iterations of $F$: $n$-th iteration of $F^*$ defines morphisms \\ $F_x^{*n}: \mathcal{F}_{F^n(x)} \to \mathcal{F}_x$; for $x$ fixed under $F^n$, $F_x^{*n}$ is an endomorphism and $$\sum_{x \in X^{F^n}} Tr(F_x^{*n}, \mathcal{F}_x)=\sum_i (-1)^i Tr(F^{*n}, H_c^i(X, \mathcal{F})) \ \ \textbf{(1.12.1)}$$

\textbf{(1.13)} Let $x_0 \in |X|$, $Z$ be the orbit corresponding to $F$ in $|X|$ and $x \in Z$. The orbit $Z$ has $\deg(x_0)$ elements (1.4). We denote by $F_{x_0}^*$ the endomorphism $F_x^{*\deg(x_0)}$ of $\mathcal{F}_x$ and let $$\det(1-F_{x_0}^*t, \mathcal{F}_0)=\det(1-F_{x_0}^*t, \mathcal{F}_x).$$ Because of the local isomorphism\footnote{The Frobenius $F$ is a local isomorphism and hence an isomorphism on stalks.}, $(\mathcal{F}_x, F_{x_0}^*)$ does not depend on the choice of $X$. This justifies omitting $x$ in the notation. We will use a similar notation for other functions of $(\mathcal{F}_x, F_{x_0}^*)$.

\textbf{(1.14)} Define $Z(X_0, \mathcal{F}_0, t) \in \mathbb{Q}_l[[t]]$ by the product $$Z(X_0, \mathcal{F}_0, t)=\prod_{x \in |X_0|} \det(1-F_x^*t^{\deg(x)}, \mathcal{F}_0)^{-1} \ \ \textbf{(1.14.1)}$$ For the constant sheaf $\mathbb{Q}_l$ we recover (1.1.2). According to (1.5.3), the logarithmic derivative of $Z$ is $$t \frac{d}{dt}\log Z(X_0, \mathcal{F}_0, t) \overset{def} = \frac{t \frac{d}{dt}Z(X_0, \mathcal{F}_0, t)}{Z(X_0, \mathcal{F}_0, t)} = \sum_n \sum_{x \in X^{F^n}=X_0(\mathbb{F}_{q^n})} Tr(F_x^{*n}, \mathcal{F}_0)t^n \ \ \textbf{(1.14.2)}$$ Substituting (1.12.1) in (1.14.2) we find by a calculation similar to (1.5) the generalization of (1.5.4) $$Z(X_0, \mathcal{F}_0, t)=\prod_i \det(1-F^*t, H_c^i(X, \mathcal{F}))^{(-1)^{i+1}} \ \ \textbf{(1.14.3)}$$ This formula is an identity in $\mathbb{Q}_l[[t]]$.

\textbf{(1.15)} It is sometimes convenient to use Galois language instead of the geometric one. Here is the dictionary. 

If $\mathbb{\bar F}_q^1$ and  $\mathbb{\bar F}_q^2$ are two algebraic closures of $\mathbb{F}_q$, $(X_0, \mathcal{F}_0)$ over $\mathbb{F}_q$ defines by extension of scalars $(X_1, \mathcal{F}_1)$ over $\mathbb{\bar F}_q^1$ and $(X_2, \mathcal{F}_2)$ over $\mathbb{\bar F}_q^2$. All $\mathbb{F}_q$-isomorphisms $\sigma: \mathbb{\bar F}_q^1 \overset{\sim} \to \mathbb{\bar F}_q^2$ define isomorphisms $$H_c^*(X_1, \mathcal{F}_1) \overset{\sim} \to H_c^*(X_2, \mathcal{F}_2).$$ In particular, for $\mathbb{\bar F}_q^1 =\mathbb{\bar F}_q^2$ (denote by $\mathbb{\bar F}_q$), we find that $Gal(\mathbb{\bar F}_q / \mathbb{F}_q)$ acts on $H_c^*(X, \mathcal{F})$ (action by transport of structure \footnote{This phrase is commonly used to state the principle that any isomorphism $Y_1 \to Y_2$ extends canonically to an isomorphism of objects constructed from $Y_1$ and $Y_2$ (cohomology groups, sheaves, etc). When $Y_1=Y_2$, automorphisms also extend.}). Let $\varphi \in Gal(\mathbb{\bar F}_q / \mathbb{F}_q)$ be the substitution of Frobenius. We verify that $$F^*=\varphi^{-1} \ \ (in \ End(H_c^*(X, \mathcal{F}))).$$ This leads to the definition of the \textit{geometric Frobenius} $F \in Gal(\mathbb{\bar F}_q / \mathbb{F}_q)$ as $\varphi^{-1}$. We have $$F^*=F \ \ \textbf{(1.15.1)}$$ Let $x$ be a geometric point of $X_0$, localized to $x_0 \in |X_0|$. By transport of structure, the group $Gal(k(x)/k(x_0))$ acts on the stalk $(\mathcal{F}_0)_x$ of $\mathcal{F}_0$ at $x$; in particular, we have a geometric Frobenius relative to $k(x_0)$: $F_{x_0} \in Gal(k(x)/k(x_0))$. For $x$ defined by a closed point, still denoted by $x$, in $X$ we have $\mathcal{F}_x=(\mathcal{F}_0)_x$ and $$F_{x_0}^* \overset{def} = F_x^{*\deg(x_0)}=F_{x_0} \ \  (in \ End(\mathcal{F}_x)) \ \ \textbf{(1.15.2)}$$ In the Galois notation, (1.14.3) looks like $$\prod_{x \in |X_0|} \det(1-F_xt^{\deg(x)}, \mathcal{F}_0)^{-1}=\prod_i \det(1-Ft, H_c^i(X, \mathcal{F}))^{(-1)^{i+1}}.$$

\section{The theory of Grothendieck: Poincare duality}
\textbf{(2.1)} To explain the relationship between the roots of unity and orientations I will first repeat the two classical cases in a wacky language.

a) \textit{Differentiable manifolds.} - Let $X$ be a differentiable manifold purely of dimension $n$. The orientation sheaf $\mathbb{Z}'$ on $X$ is the sheaf locally isomorphic to the constant sheaf $\mathbb{Z}$, whose invertible sections on an open $U$ in $X$ correspond to the orientations of $U$. An \textit{orientation} of $X$ is an isomorphism of $\mathbb{Z}'$ with the constant sheaf $\mathbb{Z}$. The \textit{fundamental class} of $X$ is a morphism $Tr: H_c^n(X, \mathbb{Z}') \to \mathbb{Z}$; if $X$ is orientable, it is identified with a morphism $Tr: H_c^n(X, \mathbb{Z}) \to \mathbb{Z}$. The Poincare duality is expressed using the fundamental class. 

b) \textit{Complex varieties.} - Let $\mathbb{C}$ be the closure of $\mathbb{R}$. A smooth complex algebraic variety or rather the underlying differentiable variety is always orientable. To justify this it suffices to orient $\mathbb{C}$ itself. This amounts to a choice: 

\textit{a)} choosing one of the two roots of the equation $X^2=-1$; we call it $+i;$

\textit{b)} choosing an isomorphism from $\mathbb{R}/\mathbb{Z}$ to $U^1=\{z \in \mathbb{C}| \ |z|=1\}$; $+i$ is the image of $\frac{1}{4}$;

\textit{c)} choosing one of the two isomorphisms $x \to \exp(\pm 2\pi ix)$ from $\mathbb{Q}/\mathbb{Z}$ to the group of the roots of unity of $\mathbb{C}$, which extends continuously to an isomorphism from $\mathbb{R}/\mathbb{Z}$ to $U^1$.

We denote by $\mathbb{Z}(1)$ a free $\mathbb{Z}$-module of rank one whose set of generators has two elements canonically corresponding to one of the two-element sets \textit{a), b), c)}. The simplest is to take $\mathbb{Z}(1)=Ker(\exp: \mathbb{C} \to \mathbb{C}^*)$. The generator $y=\pm 2 \pi i$ corresponds to the isomorphism \textit{c)}: $x \to \exp(xy)$. Let $\mathbb{Z}(r)$ be the $r$-th tensor power of $\mathbb{Z}(1)$. If $X$ is a smooth complex algebraic variety purely of complex dimension $r$, the orientation sheaf on $X$ is the constant sheaf of value $\mathbb{Z}(r)$. 

\textbf{(2.2)} To "orient" an algebraic variety over an algebraically closed $k$ of characteristic zero, we must choose an isomorphism from $\mathbb{Q}/\mathbb{Z}$ to the group of the roots of unity of $k$. The set of such isomorphisms is the principal homogeneous space for $\mathbb{\hat Z}^*$ (no longer for $\mathbb{Z}^*$). When one is only interested in the $l$-adic cohomology, it suffices to consider the roots of unity of order a power of $l$, and to assume that the characteristic $p$ of $k$ differs from $l$. We denote by $\mathbb{Z}/l^n(1)$ the group of the roots of unity of $k$ of order dividing $l^n$. For various $n$,  $\mathbb{Z}/l^n(1)$ form a projective system with transition maps $$\sigma_{m, n}: \mathbb{Z}/l^m(1) \to \mathbb{Z}/l^n(1): x \to x^{l^{m-n}}.$$ We let $\mathbb{Z}_l=\lim proj \mathbb{Z}/l^n(1)$ and $\mathbb{Q}_l(1)=\mathbb{Z}_l(1) \otimes_{\mathbb{Z}_l} \mathbb{Q}_l$. Denote by $\mathbb{Q}_l(r)$ the $r$-th tensor power of $\mathbb{Q}_l(1)$; for $r \in \mathbb{Z}$ negative we put $\mathbb{Q}_l(r)=\mathbb{Q}_l(-r \widecheck )$. 

As a vector space over $\mathbb{Q}_l$, $\mathbb{Q}_l(1)$ is isomorphic to $\mathbb{Q}_l$. However, the automorphism group of $k$ acts non-trivially on $\mathbb{Q}_l(1)$: it acts via the character with values in $\mathbb{Z}_l^*$, which gives its action on the roots of unity. In particular, for $k=\mathbb{\bar F}_q$, the substitution of Frobenius $\varphi: x \to x^q$ acts by multiplication by $q$. 

Let $X$ be an algebraic variety purely of dimension $n$ over $k$. The \textit{orientation sheaf} of $X$ for the $l$-adic cohomology is the constant $\mathbb{Q}_l$-sheaf $\mathbb{Q}_l(n)$. The \textit{fundamental class} is a morphism $$Tr: H_c^{2n}(X, \mathbb{Q}_l(n)) \to \mathbb{Q}_l,$$ or rather even $$Tr: H_c^{2n}(X, \mathbb{Q}_l) \to \mathbb{Q}_l(-n).$$

\begin{thrm2}
For $X$ proper and smooth, purely of dimension $n$, the bilinear form $$Tr(x \cup y): H^i(X, \mathbb{Q}_l) \otimes H^{2n-i}(X, \mathbb{Q}_l) \to \mathbb{Q}_l(-n)$$ is a perfect paring (it identifies $H^i(X, \mathbb{Q}_l)$ with the dual of $H^{2n-i}(X, \mathbb{Q}_l(n))$).
\end{thrm2}

\textbf{(2.4)} Let $X_0$ be a proper and smooth algebraic variety over $\mathbb{F}_q$, purely of dimension $n$ and we obtain $X$ over $\mathbb{\bar F}_q$ from $X_0$ by extension of scalars. The morphism (2.3) is compatible with the action of $Gal(\mathbb{\bar F}_q/ \mathbb{F}_q)$. If $(\alpha_j)$ are the eigenvalues of the geometric Frobenius acting on $H^i(X, \mathbb{Q}_l)$, the eigenvalues of $F$ acting on $H^{2n-i}(X, \mathbb{Q})$ are $(q^n \alpha_j^{-1})$.

\textbf{(2.5)} Assume for simplicity that $X$ is connected. The proof of (2.4) goes as follows, once we transpose to the geometric language instead of the Galois one (see (1.15)).

a) The cup-product puts $H^i(X, \mathbb{Q}_l)$ and $H^{2n-i}(X, \mathbb{Q}_l)$ into perfect duality with values in $H^{2n}(X, \mathbb{Q}_l)$ that has dimension one.

b) The cup product commutes with the inverse image of $F^*$ by the Frobenius morphism \\ $F: X \to X$.  

c) The morphism $F$ is finite of degree $q^n$: on $H^{2n}(X, \mathbb{Q}_l)$ $F^*$ is multiplication by $q^n$. 

d) Therefore, the eigenvalues of $F^*$ satisfy the property (2.4).

\textbf{(2.6)} We let $\chi(X)=\sum_i (-1)^i \dim H^i(X, \mathbb{Q}_l)$. For $n$ odd, the form $Tr(x \cup y)$ on $H^n(X, \mathbb{Q}_l)$ is skew-symmetric; the integer $n \chi(X)$ is always even. It is easy to deduce from $(1.5.4)$ and (2.3), (2.4) that $$Z(X_0, t)=\varepsilon q^{\frac{-n \chi(X)}{2}}t^{-\chi(X)}Z(X_0, q^{-n}t^{-n})$$ for $\varepsilon=\pm 1$. If $n$ is even, let $N$ denote the multiplicity of the eigenvalue $q^{\frac{n}{2}}$ of $F^*$ acting on $H^n(X, \mathbb{Q}_l)$ (i.e the dimension of the corresponding invariant subspace). We have $$\varepsilon=\begin{cases} 1, & \mbox{if } n\mbox{ is odd} \\ (-1)^N, & \mbox{if } n\mbox{ is even} \end{cases}$$ This is the Grothendieck's formulation of the functional equation for $Z$-functions. 

\textbf{(2.7)} We will need other forms of the duality theorem. The case of curves will be enough for our purposes. If $\mathcal{F}$ is a $\mathbb{Q}_l$-sheaf on an algebraic variety $X$ over an algebraically closed $k$, we denote by $\mathcal{F}(r)$ the sheaf $\mathcal{F} \otimes \mathbb{Q}_l(r)$. This sheaf is (not canonically) isomorphic to $\mathcal{F}$.

\begin{thrm3}
Let $X$ be smooth purely of dimension $n$ and $\mathcal{F}$ be a locally constant sheaf. We denote by $\mathcal{\widecheck F}$ the dual of $\mathcal{F}$. The bilinear form $$Tr(x \cup y): H^i(X, \mathcal{F}) \otimes H_c^{2n-i}(X, \mathcal{\widecheck F}(n)) \to H_c^{2n}(X, \mathcal{F} \otimes \mathcal{ \widecheck F}(n)) \to H_c^{2n}(X, \mathbb{Q}_l(n)) \to \mathbb{Q}_l$$ is a perfect pairing. 
\end{thrm3}

\textbf{(2.9)} Assume that $X$ is connected and that $x$ is a closed point of $X$. The functor $\mathcal{F} \to \mathcal{F}_x$ is an equivalence of the category of locally constant $\mathbb{Q}_l$-sheaves with that of $l$-adic representations of $\pi_1(X, x)$. Via this equivalence, $H^0(X, \mathcal{F})$ is identified with the invariants of $\pi_1(X, x)$ acting on $\mathcal{F}_x$: $$H^0(X, \mathcal{F}) \overset{\sim} \to \mathcal{F}_x^{\pi_1(X, x)}.$$ According to (2.8), for $X$ smooth and connected of dimension $n$ we have $$H_c^{2n}(X, \mathcal{F})=H^0(X, \mathcal{\widecheck F}(n) \widecheck )=((\mathcal{\widecheck F}_x(n))^{\pi_1(X, x)} \widecheck ).$$ The duality exchanges invariants (the largest invariant subspaces) with coinvariants (the largest invariant quotients)\footnote{I haven't seen the terms "invariant" and "coinvariant" used in this setting. For a group $G$ acting on $M$ the "invariants" (resp. "coinvariants") are denoted $M^G$ (resp $M_G$).}. The formula takes form $$H_c^{2n}(X, \mathcal{F})=(\mathcal{F}_x)_{\pi_1(X, x)}(-n).$$ We will only use it for $n=1$. 

\begin{stat}
Let $X$ be a connected smooth curve over an algebraically closed field $k$, $x$ a closed point of $X$ and $\mathcal{F}$ a locally constant $\mathbb{Q}_l$-sheaf. We have 

(i)$H_c^0(X, \mathcal{F})=0$ if $X$ is affine.

(ii)$H_c^{2}(X, \mathcal{F})=(\mathcal{F}_x)_{\pi_1(X, x)}(-1).$
\end{stat} Assertion (i) simply states that $\mathcal{F}$ does not have sections with finite support.

\textbf{(2.11)} Let $X$ be a connected smooth projective curve over an algebraically closed $k$, $U$ an open set in $X$, the complement of the finite set $S$ of closed points of $X$, $j$ the inclusion $U\xhookrightarrow{} X$ and $\mathcal{F}$ a locally constant $\mathbb{Q}_l$-sheaf on $U$. Let $j_*\mathcal{F}$ be the constructible $\mathbb{Q}_l$-sheaf - the direct image of $\mathcal{F}$. Its stalk at $x \in S$ has rank less or equal to the rank of the stalk at a general point; it is the space of invariants of the local monodromy group. 

\begin{thrm4}
The bilinear form $$Tr(x \cup y): H^i(X, j_*\mathcal{F}) \otimes H^{2-i}(X, j_*\mathcal{\widecheck F}(1)) \to H^{2}(X, j_*\mathcal{F} \otimes j_*\mathcal{ \widecheck F}(1)) \to$$ $$ \to H^{2}(X, j_*(\mathcal{F} \otimes \mathcal{ \widecheck F})(1)) \to H_c^{2}(X, \mathbb{Q}_l(1)) \to \mathbb{Q}_l$$ is a perfect pairing.
\end{thrm4}

\textbf{(2.13)} It will be convenient to have $\mathbb{Q}_l$-sheaves $\mathbb{Q}_l(r)$ on any scheme $X$ where $l$ is invertible. The point is to define $\mathbb{Z}/l^n(1)$. By definition, $\mathbb{Z}/l^n(1)$ is the etale sheaf of the $l^n$-th roots of unity.

\textbf{(2.14)} \textit{Bibliographical notes on paragraphs 1 and 2.} 

A) All the important results in etale cohomology are first proved for torsion sheaves. The extension to $\mathbb{Q}_l$-sheaves is done by passing to formal limits. In what follows, for each theorem mentioned I will not refer to the reference where it is proved, but to the reference where a similar statement for a torsion sheaf is proved. 

B) With the exception of the Lefschetz formula and (2.12), results in etale cohomology used in this article are all proved in SGA 4. For those already stated, the references are: definition of $H^i$: VII; definition of $H_c^I$: XVII 5.1; finiteness theorem: XIV 1, completed in XVII 5.3; cohomological dimension: X; Poincare duality: XVIII.

C) The relation between the various Frobeniuses ((1.4), (1.11), (1.15)) is explained in detail in SGA 5, XV, par. 1, 2.

D) The cohomological interpretation of the $Z$-functions is clearly exposed in [1]; however, Lefschetz formula (1.12) for $X$ a smooth projective curve is used, but not proved. For the proof, one has to consult SGA 5. 

E) The form $(2.12)$ of the Poincare duality theorem follows from the general result SGA 4, XVIII (3.2.5) (for $S=Spec(k), \ X=X, \ K=j_*\mathcal{F}, \ L=\mathbb{Q}_l$) by a local calculation that is not difficult. The statement will be explicitly included in the final version of SGA 5. For the case interesting to us (tame ramification of $\mathcal{F}$), we could obtain it transcendentially\footnote{Transcendential algebraic geometry deals with varieties defined over $\mathbb{C}$ and concentrates on their structure of holomorphic manifolds, that allows one to use powerful techniques of topology, analysis, differential equations, etc.} by lifting $X$ and $\mathcal{F}$ to characteristic $0$. 

\section{The main lemma (La majoration fondamentale)}
The result of this paragraph was catalyzed by reading the lecture of Rankin [3]\footnote{For a slightly different exposition of the main lemma using the equivalence of (1.9) b) consult Milne.}. 

\textbf{(3.1)} Let $U_0$ be a curve on $\mathbb{F}_q$, complement in $\mathbb{P}^1$ to a finite set of closed points, $U$ be the curve over $\mathbb{\bar F}_q$ corresponding to it, $u$ a closed point of $U$, $\mathcal{F}_0$ a locally constant sheaf on $U_0$ and $\mathcal{F}$ its inverse image on $U$. 

Let $\beta \in \mathbb{Q}$. We say that $\mathcal{F}_0$ is of \textit{weight} $\beta$ if for all $x \in |U_0|$, the eigenvalues of $F_x$ acting on $\mathcal{F}_0$ (1.13) are algebraic numbers all of which complex conjugates are of absolute value $q_x^{\beta/2}$. For example, $\mathbb{Q}_l(r)$ is of weight $-2r$.

\begin{thrm5}
Let's make the following hypotheses:

(i) $\mathcal{F}_0$ is equipped with a bilinear skew-symmetric nondegenerate form $$\psi: \mathcal{F}_0 \otimes \mathcal{F}_0 \to \mathbb{Q}_l(-\beta) \ \ (\beta \in \mathbb{Z}).$$

(ii) The image of $\pi_1(U, u)$ in $GL(\mathcal{F}_u)$ is an open subgroup of the symplectic group $Sp(\mathcal{F}_u, \psi_u)$.

(iii) For all $x \in |U_0|$, the polynomial $\det(1-F_x t, \mathcal{F}_0)$ has rational coefficients. 

Then $\mathcal{F}$ is of weight $\beta$.

\end{thrm5}
We may and do assume that $U$ is affine and that $\mathcal{F} \neq 0$. 

\begin{lemma1}
Let $2k$ be an even integer and denote by $\overset{2k} \otimes \mathcal{F}_0$ the $2k$-th tensor power of $\mathcal{F}_0$. For $x \in |U_0|$ the logarithmic derivative $$t \frac{d}{dt}\log(\det(1-F_x t^{\deg(x)}, \overset{2k} \otimes \mathcal{F}_0)^{-1})$$ is a formal series with positive rational coefficients.  
\end{lemma1}

Hypothesis (iii) ensures that for all $n$  $Tr(F_x^n, \mathcal{F}_0) \in \mathbb{Q}$. The number $$Tr(F_x^n, \overset{2k} \otimes \mathcal{F}_0)=Tr(F_x^n, \mathcal{F}_0)^{2k}$$ is a positive rational and we apply (1.5.3). 

\begin{lemma2}
The local factors $\det(1-F_x t^{\deg(x)}, \overset{2k} \otimes \mathcal{F}_0)^{-1}$ are formal series with positive rational coefficients.
\end{lemma2}

The formal series $\log(\det(1-F_x t^{\deg(x)}, \overset{2k} \otimes \mathcal{F}_0)^{-1})$ has constant term zero, from (3.3) all the coefficients are $\geq 0$; the coefficients of the exponentiation are therefore also positive.

\begin{lemma3}
Let $f_i=\sum_n a_{i, n}t^n$ be a sequence of formal series with positive real coefficients. We assume that the order of $f_i-1$ tends to infinity with $i$; and we denote $f=\prod_i f_i$. Then the radius of convergence of $f_i$ is greater or equal to that of $f$. 
\end{lemma3}

If $f=\sum_n a_n t^n$, we have $a_{i, n} \leq a_n$.

\begin{lemma4}
Under the assumptions of (3.5), if $f$ and the $f_i$ are Taylor expansions of meromorphic functions, then $$inf \{ |z| \ | f(z)=\infty \} \leq inf \{ |z| \ | f_i(z)=\infty \}$$
\end{lemma4}

Indeed, those numbers are the radii of convergence. 

\textbf{(3.7)} For each partition $P$ of $[1, 2k]$ \footnote{of $\{1, \cdots, 2k\}$} into the two element sets $\{ i_{\alpha}, j_{\alpha} \} (i_{\alpha} \leq j_{\alpha})$, we define $$\psi_P: \overset{2k} \otimes \mathcal{F}_0 \to \mathbb{Q}_l(-k\beta): x_1 \otimes \cdots \otimes x_{2k} \to \prod_{\alpha} \psi(x_{i_\alpha}, x_{j_{\alpha}}).$$ Let $x$ be a closed point of $X$ \footnote{Let $u$ be a closed point of $U$.}. Hypothesis (ii) ensures that the coivariants of $\pi_1(U, u)$ on $\overset{2k} \otimes \mathcal{F}_u$ are the coinvariants on $\overset{2k} \otimes \mathcal{F}_u$ of the entire symplectic group ($\pi_1$ is Zariski-dense in $Sp$). Let $\mathcal{P}$ be the set of partitions $P$. From H.Weil (\textit{The classical groups}, Princeton University Press, chap. VI, par. 1), for an appropriate $\mathcal{P}' \subset \mathcal{P}$, depending on $\dim(\mathcal{F}_u)$, $\psi_P$ (for $P \in \mathcal{P}'$) defines an isomorphism $$(\overset{2k} \otimes \mathcal{F}_u)_{\pi_1}=(\overset{2k} \otimes \mathcal{F}_u)_{Sp} \overset{\sim} \to \mathbb{Q}_l(-k\beta)^{\mathcal{P}'}.$$ Let $N$ be the number of elements in $\mathcal{P}'$. According to (2.10) the formula above gives $$H_c^2(U, \overset{2k} \otimes \mathcal{F}) \simeq \mathbb{Q}_l(-k\beta-1)^N.$$ Since $H_c^0(U, \overset{2k} \otimes \mathcal{F})=0$ \footnote{Again, see (2.10).}, the formula (1.14.3) reduces to $$Z(U_0, \overset{2k} \otimes \mathcal{F}_0, t)=\frac{\det(1-F^*t, H^1(U, \overset{2k} \otimes \mathcal{F}))}{(1-q^{k\beta+1}t)^N}.$$ This $Z$-function is therefore the Taylor series expansion of a rational function having only one pole at $t=1/q^{k\beta+1}$. We will only use the fact that the poles are of modulus $t=1/q^{k\beta+1}$ in $\mathbb{C}$. This could be concluded from the general arguments about reductive groups. If $\alpha$ is an eigenvalue of $F_x$ on $\mathcal{F}_0$, then $\alpha^{2k}$ is an eigenvalue of $F_x$ on $\overset{2k} \otimes \mathcal{F}_0$. We now let $\alpha$ be any complex conjugate of \footnote{the original} $\alpha$. The inverse power $1/\alpha^{2k/\deg(x)}$ is a pole of $\det(1-F_x t^{\deg(x)}, \overset{2k} \otimes \mathcal{F})^{-1}$. According to (3.4) and (3.6) we have $$|1/q^{k\beta+1}| \leq |1/\alpha^{2k/\deg(x)}|,$$ or $$|\alpha| \leq q_x^{\frac{\beta}{2}+\frac{1}{2k}}.$$ Letting $k$ go to infinity we find that $$|\alpha| \leq q_x^{\frac{\beta}{2}}.$$ On the other hand, the existence of $\psi$ ensures that $q_x^\beta \alpha^{-1}$ is an eigenvalue, so we have the inequality $$|q_x^\beta \alpha^{-1}| \leq q_x^{\beta/2},$$ or $$q_x^{\beta/2} \leq |\alpha|.$$ This completes the proof.

\begin{cor}
Let $\alpha$ be an eigenvalue of $F^*$ acting on $H_c^1(U, \mathcal{F})$. Then $\alpha$ is an algebraic number all of which complex conjugates satisfy $$|\alpha| \leq q^{\frac{\beta+1}{2}+\frac{1}{2}}.$$

\end{cor}

The formula (1.14.3) for $\mathcal{F}_0$ reduces to $$Z(U_0, \mathcal{F}, t)=\det(1-F^*t, H_c^1(U, \mathcal{F})).$$ The left hand side is a formal series with rational coefficients, given its representation as a product and hypothesis (iii). The right hand side is therefore a polynomial with rational coefficients, $1/\alpha$ is its root. This already implies that $\alpha$ is algebraic. To complete the proof, it suffices to show that the infinite product that defines $Z(U_0, \mathcal{F}_0, t)$ converges absolutely (thus it is nonzero) for $|t| <  q^{\frac{-\beta}{2}-1}$. 

Let $N$ be the rank of $\mathcal{F}$ and let $$\det(1-F_xt, \mathcal{F})=\prod_{i=1}^N (1-\alpha_{i, x}t).$$ According to (3.2), $|\alpha_{i, x}|=q_x^{\beta/2}$. The convergence of the infinite product for Z follows from that of the series $$\sum_{i, x} |\alpha_{i, x} t^{\deg(x)}|.$$ For $|t|=q^{\frac{-\beta}{2}-1-\varepsilon}\ (\varepsilon >0)$ we have $$\sum_{i, x} |\alpha_{i, x} t^{\deg(x)}|=N\sum_x q_x^{-1-\varepsilon}.$$ On the affine line there are $q^n$ points with coordinate in $\mathbb{F}_{q^n}$, so there are at most $q^n$ closed points of degree $n$. So we have $$\sum_x q_x^{-1-\varepsilon} \leq \sum_n q^nq^{n(-1-\varepsilon)}=\sum_n q^{-n\varepsilon} < \infty,$$ which completes the proof.

\begin{cor1}
Let $j_0$ be the inclusion of $U_0$ in $\mathbb{P}_{\mathbb{F}_q}^1$, $j$ that of $U$ into $\mathbb{P}^1$ and $\alpha$ an eigenvalue of $F^*$ acting on $H^1(\mathcal{P}^1, j_*\mathcal{F})$. Then $\alpha$ is an algebraic number all of which complex conjugates satisfy $$q^{\frac{\beta+1}{2}-\frac{1}{2}} \leq |\alpha| \leq q^{\frac{\beta+1}{2}+\frac{1}{2}}.$$
\end{cor1}

The segment of the long exact sequence in cohomology defined by the short exact sequence  $$0 \to j_! \mathcal{F} \to j_*\mathcal{F} \to j_*\mathcal{F}/j_! \mathcal{F} \to 0$$ ($j_!$ is the extension by 0) is $$H_c^1(U, \mathcal{F}) \to H^1( \mathbb{P}, j_* \mathcal{F}) \to 0.$$ So the eigenvalue $\alpha$ already appears in $H_c^1(U, \mathcal{F})$ \footnote{We will repeatedly use the following facts: 

A) If $H^i(X') \to H^i(X)$ is a surjection, then (1.7) for $H^i(X')$ $\Rightarrow$ (1.7) for $H^i(X)$.  

B) If $H^i(X) \xhookrightarrow{} H^i(X'')$ is an embedding, then (1.7) for $H^i(X'')$ $\Rightarrow$ (1.7) for $H^i(X)$.} and so by (3.8) we have: $$|\alpha| \leq q^{\frac{\beta+1}{2}+\frac{1}{2}}.$$ The Poincare duality (2.12) implies that $q^{\beta+1}\alpha^{-1}$ is an eigenvalue, so we have the inequality $$|q^{\beta+1}\alpha^{-1}| \leq q^{\frac{\beta+1}{2}+\frac{1}{2}}$$ and the corollary is proved.

\section{Lefschetz theory: local theory}

\textbf{(4.1)} On $\mathbb{C}$ Lefshietz local results are as follows. Let $D=\{z| \ |z| < 1\}$ be the unit disk, $D^*=D-\{0\}$ and $f:X \to D$ be a morphism of analytic spaces. We assume that

a) $X$ is nonsingular and purely of dimension $n+1$;

b) $f$ is proper;

c) $f$ is smooth outside the point $x$ of the special fiber $X_0=f^{-1}(0)$;

d) At $x$ $f$ has a nondegenerate double point.

Let $t \neq 0$ in $D$ and $X_t=f^{-1}(t)$ "the" general fiber. To the previous data we associate:

$\alpha$) Specialization morphisms $sp: H^i(X_0, \mathbb{Z}) \to H^i(X_t, \mathbb{Z})$: $X_0$ is a deformation retract of $X$ and $sp$ is the composition arrow $$H^i(X_0, \mathbb{Z}) \overset{\sim} \leftarrow H^i(X, \mathbb{Z})\to H^i(X_t, \mathbb{Z})$$

$\beta$) The monodromy transformations $T: H^i(X_t, \mathbb{Z}) \to H^i(X_t, \mathbb{Z})$, which describe the effect on the singular cycles of $X_t$ of "rotating $t$ around $0$". This is even an action on $H^i(X_t, \mathbb{Z})$, the stalk at $t$ of the local system $R^i f_* \mathbb{Z}|D^*$, of the positive generator of $\pi_1(D^*, t)$. 

Lefschetz theory describes $\alpha$) and $\beta$) in terms of the \textit{vanishing cycle}\footnote{We can define $\delta$ to be the unique (up to sign) generator of $H^n(X_0, \mathbb{Z})^{\perp} \subset H^n (X_t, \mathbb{Z})$ (under the pairing induced by Poincare duality).} $\delta \in H^n (X_t, \mathbb{Z})$. This cycle is well-defined up to sign. For $i \neq n, \ n+1$ we have $$H^i(X_0, \mathbb{Z}) \overset{\sim} \to H^i(X_t, \mathbb{Z}) \ (i \neq n, \ n+1).$$  For $i=n, \ n+1$ we have an exact sequence $$0 \to H^n(X_0, \mathbb{Z}) \to H^n(X_t, \mathbb{Z}) \overset{x \to (x, \delta)} \to \mathbb{Z} \to H^{n+1}(X_0, \mathbb{Z}) \to H^{n+1}(X_t, \mathbb{Z}) \to 0.$$ For $i \neq n$, the monodromy $T$ is the identity. For $i=n$ we have $$Tx=x \pm (x, \delta)\delta.$$ 
The values of $\pm$, $T\delta$ and $(\delta, \delta)$ are as follows: 
\begin{center}
\begin{tabular}{ c c c c c}
 $n \mod 4$ & 0 & 1 & 2 & 3 \\
 $Tx=x\pm(x, \delta)\delta$ & - & - & + & + \\  
$(\delta, \delta)$ & 2 & 0 & -2 & 0 \\ $T\delta$ & $-\delta$ & $\delta$ & $-\delta$ & $\delta$    
\end{tabular}
\end{center}
The monodromy transformation preserves the intersection form $Tr(x \cup y)$ on $H^n(X_t, \mathbb{Z})$. For $n$ odd, it is the symplectic transvection. For $n$ even, it is the symmetric orthogonal. 

\textbf{(4.2)} There is an analog of (4.1) in abstract algebraic geometry. The disk $D$ is replaced by the spectrum of a henselian discrete valuation ring $A$ with an algebraically closed residue field. Let $S$ be the spectrum, $\eta$ its generic point (spectrum of the field of fractions of $A$), $s$ the closed point (spectrum of the residue field). The role of $t$ is played by the geometric generic point $\bar \eta$ (spectrum of the closure of the field of fractions of $A$). 

Let $f: X \to S$ be a proper morphism, with $X$ regular purely of dimension $n+1$. We assume that $f$ is smooth except for an ordinary double point $x$ of the special fiber $X_s$. Let $l$ be a prime number different from the residual characteristic\footnote{Characteristic of the residue field of $A$.}  $p$ of $S$. Denoting by $X_{\bar \eta}$ the generic geometric fiber, we have the specialization morphism $$sp: H^i(X_s, \mathbb{Q}_l) \overset{\sim} \leftarrow H^i(X, \mathbb{Q}_l) \to H^i(X_{\bar \eta}, \mathbb{Q}_l) \ \ \textbf{(4.2.1)}$$ The role of $T$ is played by the action of the inertia group $I=Gal(\bar \eta/ \eta)$ on $H^i(X_{\bar \eta}, \mathbb{Q}_l)$ by transport of structure (see (1.15)): $$I=Gal(\bar \eta/\eta) \to GL(H^i(X_{\bar \eta}, \mathbb{Q}_l)) \ \ \textbf{(4.2.2)}$$ The data (4.2.1), (4.2.2) fully determines the sheaf $R^i f_*\mathbb{Q}_l$ on $S$. 

\textbf{(4.3)} Let $n=2m$ for $n$ even and $n=2m+1$ for $n$ odd. (4.2.1) and (4.2.2) can still be described in terms of the vanishing cycle $$\delta \in H^n(X_{\bar \eta}, \mathbb{Q}_l)(m) \ \ \textbf{(4.3.1)}$$ This cycle is well-defined up to sign. 

For $i \neq n, \ n+1$ we have $$H^i(X_s, \mathbb{Q}_l) \overset{\sim} \to H^i(X_{\bar \eta}, \mathbb{Q}_l) \ (i \neq n, \ n+1) \ \ \textbf{(4.3.2)}$$  For $i=n, \ n+1$ we have an exact sequence $$0 \to H^n(X_s, \mathbb{Q}_l) \to H^n(X_{\bar \eta}, \mathbb{Q}_l) \overset{x \to Tr(x \cup \delta)} \to \mathbb{Q}_l(m-n) \to H^{n+1}(X_s, \mathbb{Q}_l) \to H^{n+1}(X_{\bar \eta}, \mathbb{Q}_l) \to 0 \ \ \textbf{(4.3.3)}$$

The action (4.2.2) (local monodromy) is trivial for $i \neq n$. For $i=n$, it is described as follows.

A) \textit{n odd}. - We have a canonical homomorphism $$t_l: I \to \mathbb{Z}_l(1),$$ and the action of $\sigma \in I$ is $$x \to x\pm t_l(\sigma)(x, \delta)\delta.$$

B) \textit{n even}. - We will not need this case. Let's just say that if $p \neq 2$, there exists a unique character of order two $$\varepsilon: I \to \{\pm1\},$$ and we have 
\begin{center}
\begin{tabular}{ c c c }
 $\sigma x=x$ & if & $\varepsilon(\sigma)=1$ \\ $\sigma x=x\pm(x, \delta)\delta$ & if & $\varepsilon(\sigma)=-1$    
\end{tabular}
\end{center}
The signs $\pm$ in A) and B) are the same as in (4.1).

\textbf{(4.4)} These results imply the following information\footnote{Derivation of some of the results can be found in Milne.} about $R^if_*\mathbb{Q}_l$.

a) If $\delta \neq 0$:

\ \ \ \ 1) For $i \neq n$ the sheaf $R^if_*\mathbb{Q}_l$ is constant.

\ \ \ \ 2) Let $j$ be the inclusion of $\eta$ in $S$. We have $$R^if_*\mathbb{Q}_l=j_* j^* R^if_*\mathbb{Q}_l.$$

b) If $\delta=0$: (This is an exceptional case. Since $(\delta, \delta)=\pm 2$ for $n$ even, it can only happen for $n$ odd.)

\ \ \ \ 1) For $i \neq n+1$ the sheaf $R^if_*\mathbb{Q}_l$ is constant.

\ \ \ \ 2) Let $\mathbb{Q}_l(m-n)_s$ be the sheaf $\mathbb{Q}_l(m-n)$ on $\{s\}$, extended by zero on $S$. Then we have an exact sequence $$0 \to \mathbb{Q}_l(m-n)_s \to R^{n+1}f_* \mathbb{Q}_l \to j_* j^*R^{n+1}f_* \mathbb{Q}_l \to 0,$$ where $j_* j^*R^{n+1}f_* \mathbb{Q}_l$ is a constant sheaf. 

\section{Lefschetz theory: global theory}
\textbf{(5.1)} On $\mathbb{C}$ the results of Lefshietz are as follows. Let $\mathbb{P}$ be a projective space of dimension $\geq 1$ and $ \mathbb{\widecheck P}$ the dual projective space; its points parameterize the hyperplanes of $\mathbb{P}$ and we denote by $H_t$ the hyperplane defined by $t \in \mathbb{\widecheck P}$. If $A$ is  a linear subspace of codimension 2 in $\mathbb{P}$, the hyperplanes containing $A$ are parameterized by points of the line $D \subset \mathbb{\widecheck P}$, the \textit{dual} of $A$. These hyperplanes $(H_t)_{t \in D}$ form the \textit{pencil with axis} $A$. 

Let $X \subset \mathbb{P}$ be a connected nonsingular projective variety of dimension $n+1$. Let $\tilde X \subset X \times D$ be the set of pairs $(x, t)$ such that $x \in H_t$. Projections to the first and second coordinates form a diagram\footnote{In the original paper the letter $f$ is absent.} $$\begin{tikzcd}
X \arrow[r, leftarrow, "\pi"]
& \tilde X \arrow[d, "f"]\\
& D
\end{tikzcd} \ \ \textbf{(5.1.1)}$$ The fiber of $f$ at $t \in D$ is the hyperplane section $X_t=X \cap H_t$ of $X$.

Fix $X$ and take a general enough $A$. Then:

A) $A$ is transverse to $X$ and $\tilde X$ is the blowing up of $X$ along $A \cap X$. In particular, $\tilde X$ is nonsingular.

B) There is a finite subset $S$ of $D$ and for each $s \in S$ a point $x_s \in X_s$ such that $f$ is smooth outside $x_s$.

C) The $x_s$ are critical nondegenerate points of $f$. 

Therefore, for each $s \in S$ local Lefschetz theory (4.1) applies to a small disk $D_s$ around $s$ and $f^{-1}(D_s)$. 

\textbf{(5.2)} Let $U=D-S$. Let $u \in U$ and choose disjoint loops $(\gamma_s)_{s \in S}$ starting from $u$, with $\gamma_s$ turning once around $s$: 
\begin{figure}[h!]
\begin{center}
\includegraphics[scale=0.5]{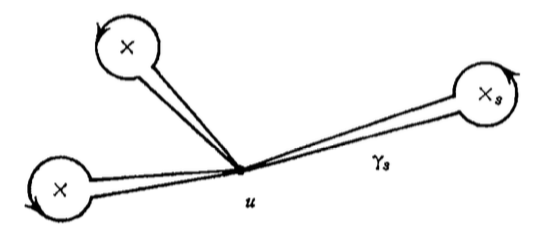}
\end{center}
\end{figure}

These loops generate the fundamental group $\pi_1(U, u)$. This group acts on $H^i(X_u, \mathbb{Z})$, the stalk at $u$ of the local system $R^if_*\mathbb{Z}| U$. According to the local theory (4.1), to each $s \in S$ corresponds a vanishing cycle $\delta_s \in H^n(X_u, \mathbb{Z})$; these cycles depend on the choice of $\gamma_s$. For $i \neq n$, the action of $\pi_1(U, u)$ on $H^i(X_u, \mathbb{Z})$ is trivial. For $i=n$ we have $$\gamma_s x=x \pm (x, \delta_s)\delta_s \ \ \textbf{(5.2.1)}$$

Let $E$ be the subspace of $H^n(X_u, \mathbb{Q})$ generated by the $\delta_s$ (\textit{vanishing} part of the cohomology). 

\begin{fact1}
$E$ is stable under the action of the monodromy group $\pi_1(U, u)$. The orthogonal $E^{\perp}$ of $E$ (for the intersection form $Tr(x \cup y)$) is the space of the invariants of the monodromy in $H^n(X_u, \mathbb{Q})$. 
\end{fact1}
The $\gamma_s$ generate the monodromy group, so this is clear from (5.2.1).

\begin{thrm6}
The vanishing cycles $\pm \delta_s$ are conjugate (up to sign)\footnote{That is, given $s, \ s' \in S$, there exists a $\sigma \in \pi_1$ such that $\sigma \delta_s=\pm \delta_s$.} under the action of $\pi_1(U, u)$. 
\end{thrm6}

Let $\widecheck X \subset \mathbb{\widecheck P}$ be the dual variety of $X$; it is the set of $t \in \mathbb{\widecheck P}$ such that $H_t$ is tangent to $X$, i.e. such that $X_t$ is singular or $X \subset H_t$. The variety $\widecheck X$ is irreducible. Let $Y \subset X \times \mathbb{\widecheck P}$ be the space of pairs $(x, t)$ such that $x \in H_t$. We have a diagram $$\begin{tikzcd}
X \arrow[r, leftarrow]
& Y \arrow[d, "g"]\\
& \mathbb{\widecheck P}
\end{tikzcd}$$ The fiber of $g$ at $t \in \mathbb{\widecheck P}$ is the hyperplane section $X_t= X \cap H_t$ of $X$ and $g$ is smooth on the complement of the inverse image of $\widecheck X$. 

We retrieve the situation of (5.1) by replacing $\mathbb{\widecheck P}$ by the line $D \subset \mathbb{\widecheck P}$ and $Y$ by $g^{-1}(D)$. We have $S=D \cap \widecheck X$. According to a theorem of Lefschetz, for $D$ general enough, the map $$\pi_1(D-S, u) \to \pi_1(\mathbb{\widecheck P}-\widecheck X, u)$$ is surjective. It suffices to show that $\pm \delta_s$ are conjugate under $\pi_1(\mathbb{\widecheck P}- \widecheck X)$. 

For $x$ in the smooth locus of codimension $1$ of $\widecheck X$, let $ch$ be the path from $t$ to $x$ in $\mathbb{\widecheck P}- \widecheck X$ and $\gamma_x$ the loop that follows $ch$ until the neighborhood of $\widecheck X$, turns once around $\widecheck X$ and then returns to $t$ by $ch$. The loops $\gamma_x$ (for various $ch$) are mutually conjugate. Since $\widecheck X$ is irreducible, two points in the smooth locus of $\widecheck X$ can always be joined, in $\widecheck X$, by a path that does not leave the smooth locus. It follows that the conjugation class of $\gamma_x$ does not depend on $x$. In particular, $\gamma_s$ are mutually conjugate. We see from (5.2.1) that this implies the conjugacy of $\pm \delta_s$. 

\begin{cor2}
The action of $\pi_1(U, u)$ on $E/(E \cap E^{\perp})$ is absolutely irreducible\footnote{The action on the $k$-vector space $V$ is called absolutely irreducible if the corresponding action on $V \otimes_k \bar k$ is irreducible.}.
\end{cor2}

Let $F \subset E \otimes \mathbb{C}$ be the subspace stable under the monodromy. If $F \not\subset (E \cap E^{\perp}) \otimes \mathbb{C}$, there exists an $x \in F$ and $s \in S$ such that $(x, \delta_s) \neq 0$. We then have $$\gamma_s x-x=\pm (x, \delta_s)\delta_s \in F$$ and $\delta_s \in F$. According to (5.4), all the $\delta_s$ are then in $F$ and $F=E$. This proves (5.5). 

\textbf{(5.6)} These results transpose as follows into abstract algebraic geometry. 
Let $\mathbb{P}$ be a projective space of dimension $>1$ over an algebraically closed field $k$ of characteristic $p$ and $X \subset \mathbb{P}$ a connected projective nonsingular variety of dimension $n+1$. For $A$ a linear subspace of codimension 2 we define $D$, the pencil $(H_t)_{t \in D}$, $\tilde X$ and the diagram (5.1.1) as in (5.1). We say that $(H_t)_{t \in D}$ form a \textit{Lefschetz pencil} of hyperplane sections if the following conditions are satisfied:

A) The axis $A$ is transverse to $X$. $\tilde X$ is obtained by blowing up $X$ along $A \cap X$ and is smooth. 

B) There is a finite subset $S$ of $D$ and for each $s \in S$ a point $x_s \in X_s$ such that $f$ is smooth outside $x_s$.

C) The $x_s$ are ordinary double singular points of $f$. 

For each $s \in S$ the local Lefschetz theory of par. 4 applies to the spectrum $D_s$ of the henselization of the local ring of $D$ at $s$ and to $\tilde X_{D_s}=\tilde X \times_D D_s$. 

\textbf{(5.7)} Let $N$ be the dimension of $\mathbb{P}$, $r$ an integer $\geq 1$ and $\iota_{(r)}$ the embedding of $\mathbb{P}$ into the projective space of dimension ${N+r \choose N} -1$, the homogeneous coordinates of which are monomials of degree $r$ in the homogeneous coordinates of $\mathbb{P}$ \footnote{Deligne describes the Veronese mapping.}. The hyperplane sections of $\iota_{(r)}(\mathbb{P})$ are the hypersurfaces of degree $r$ of $\mathbb{P}$.

For $p \neq 0$ it might happen that there is no such pencil of hyperplane sections of $X$ that is Lefschetz. However, if $r \geq 2$ and we replace the projective embedding $\iota_1: X \xhookrightarrow{} \mathbb{P}$ by $\iota_r=\iota_{(r)} \circ \iota_1$, then, in this new embedding, any general enough pencil of hypersurface sections of degree $r$ on $X$ is still Lefschetz.

\textbf{(5.8)} For the rest of this discussion, we are studying the Lefschetz pencil of hyperplane sections of $X$, \textit{excluding the case $p=2$, n even}. The case of $n$ odd will suffice for our purposes. We put $U=D-S$. Take $u \in U$ and $l$ a prime number $\neq p$. The local results of par. 4 show that $R^n f_* \mathbb{Q}_l$ is \textit{tamely} ramified at each $s \in S$. The tame fundamental group of $U$ is a quotient of the profinite completion of the corresponding transcendential\footnote{That is, taken over $\mathbb{C}$, see footnote 7.} fundamental groups (lifting to characteristic $0$ of the tame coverings and the Riemann existence theorem). The algebraic situation is therefore similar to the transcendential situation and the transfer of Lefschetz's results is done by standard arguments. In the proof of (5.4) the theorem of Lefschetz for $\pi_1$ is replaced by a theorem of Bertiny and we have to invoke Abhyankar's lemma to control the ramification of $R^*g_*\mathbb{Q}_l$ along the smooth locus of codimension one in $\widecheck X$ \footnote{See the sketch of this proof in Milne.}. 

The results are as follows\footnote{Again, proofs of some of the results appear in Milne.}:

a) \textit{If the vanishing cycles are nonzero}: 

\ \ \ \ 1) For $i \neq n$ the sheaf $R^if_*\mathbb{Q}_l$ is constant.

\ \ \ \ 2) Let $j$ be the inclusion of $U$ in $D$. We have $$R^nf_*\mathbb{Q}_l=j_* j^* R^nf_*\mathbb{Q}_l.$$

\ \ \ \ 3) Let $E \subset H^n(X_u, \mathbb{Q}_l)$ be the subspace of the cohomology generated by the vanishing cycles. This subspace is stable under $\pi_1(U, u)$ and $$E^{\perp}=H^n(X_u, \mathbb{Q}_l)^{\pi(U, u)}.$$ The representation of $\pi(U, u)$ on $E/(E \cap E^{\perp})$ is absolutely irreducible and the image of $\pi_1$ in $GL(E/(E \cap E^{\perp}))$ is generated (topologically\footnote{They generate a dense subgroup.}) by the $x \to x\pm(x, \delta_s)\delta_s \ (s \in S)$ (the $\pm$ sign is determined as in (4.1)). 

b) \textit{If the vanishing cycles are zero}: (This is an exceptional case. Since $(\delta, \delta)=\pm 2$ for $n$ even, it can only happen for $n$ odd: $n=2m+1$. Note that if one vanishing cycle is zero, they all are because of conjugacy.)

\ \ \ \ 1) For $i \neq n+1$ the sheaf $R^if_*\mathbb{Q}_l$ is constant.

\ \ \ \ 2) We have an exact sequence $$0 \to \underset{s \in S} \oplus \mathbb{Q}_l(m-n)_s \to R^{n+1}f_*\mathbb{Q}_l \to \mathcal{F} \to 0$$ with $\mathcal{F}$ constant. 

\ \ \ \ 3) $E=0$.

\textbf{(5.9)} The subspace $E \cap E^{\perp}$ of $E$ is the kernel of the restriction to $E$ of the intersection form $Tr(x \cup y)$. Therefore, this form induces a bilinear nondegenerate form $$\psi: E/(E \cap E^{\perp}) \otimes E/ (E \cap E^{\perp}) \to \mathbb{Q}_l(-n),$$ skew-symmetric for $n$ odd and symmetric for $n$ even. This form is preserved by the monodromy; for $n$ odd, therefore, the monodromy representation induces $$\rho: \pi_1(U, u) \to Sp(E/(E \cap E^{\perp}), \psi).$$

\begin{thrm7}
\footnote{See a somewhat more detailed exposition of this part in Milne.} The image of $\rho$ is open. 
\end{thrm7}

The image of $\rho$  is a compact, therefore, analytic $l$-adic\footnote{A Lie group over $\mathbb{Q}_l$.} subgroup of $Sp(E/(E \cap E^{\perp}), \psi)$. It suffices to show that its Lie algebra $\mathfrak{L}$ equals $\mathfrak{sp}(E/(E \cap E^{\perp}), \psi)$. The transcendential analog of this Lie algebra is the Lie algebra of the Zariski closure of the monodromy group. 

We deduce from (5.8) that $\mathfrak{L}$ is generated by transformations with zero square $$N_s: x \to (x, \delta_s)\delta_s \ \ (s \in S)$$ and that $E/(E \cap E^{\perp})$ is an absolutely irreducible representation of $\mathfrak{L}$. The theorem comes from the following lemma. 

\begin{lemma5}
Let $V$ be a finite dimensional vector space over the field $k$ of characteristic $0$, $\psi$ a nondegenerate skew-symmetric form on a Lie subalgebra $\mathfrak{L}$ of the Lie algebra $\mathfrak{sp}(V, \eta)$. We assume that: 

\ \ \ \ (i) $V$ is a simple representation of $\mathfrak{L}$. 

\ \ \ \ (ii) $\mathfrak{L}$ is generated by the family of endomorphisms of $V$ of the form $x \to \psi(x, \delta)\delta$. \\
Then $\mathfrak{L}=\mathfrak{sp}(V, \psi)$.
\end{lemma5}
We may and do assume that $V$, and thus $\mathfrak{L}$ are nonzero. Let $W \subset V$ be the set of $\delta \in V$ such that $N(\delta): x \to \psi(x, \delta)\delta$ is in $\mathfrak{L}$. 

a) $W$ is stable under homotheties (since $\mathfrak{L}$ is a vector subspace of $\mathfrak{gl}(V)$). 

b) If $\delta \in W$, $\exp(\lambda N(\delta))$ is an automorphism of $(V, \psi, \mathfrak{L})$, therefore, transforms $W$ to itself. If $\delta', \ \delta'' \in W$, we have $\exp(\lambda N(\delta'))\delta''=\delta''+\lambda \psi(\delta'', \delta')\delta' \in W$ \footnote{We use the nilpotence of $N(\delta)$ first to define $\exp(\lambda N(\delta))$ and then to obtain the formula.}; if $\psi(\delta', \delta'') \neq 0$, then the vector subspace spanned by $\delta'$ and $\delta''$ lies in $W$. 

c) It follows that $W$ is the union of its maximal linear subspaces $W_{\alpha}$ and that those are pairwise orthogonal.  Each $W_{\alpha}$ is therefore stable under the $N(\delta) \  (\delta \in W)$, so it is stable under $\mathfrak{L}$. By hypothesis (i), $W_{\alpha}=V$ and $\mathfrak{L}$ contains all $N(\delta)$ for $\delta \in V$.  We conclude by noting that Lie algebra $\mathfrak{sp}(V, \psi)$ is generated by the $N(\delta) \ (\delta \in V)$. 

\textit{Remark} \textbf{(5.12)} (not necessary for the exposition). - It is now easy to prove (1.6) for a hypersurface of odd dimension $n$ in $\mathbb{P}_{\mathbb{F}_q}^{n+1}$. 

Let $X_0$ be such a hypersurface and $\bar X_0$ the hypersurface over $ \mathbb{\bar F}_q$, which is obtained by extension of scalars. We have $$H^i(\bar X_0, \mathbb{Q}_l)=\mathbb{Q}_l(-i) \ \ (0 \leq i \leq n);$$ $H^i(\bar X_0, \mathbb{Q}_l(i))$ is generated by the $i$-th cup power of $\eta$, the cohomology class $c_1(\mathcal{O}(1))$ of a hyperplane section. Therefore, we have $$Z(X_0, t)=\det(1-F^*t, H^n(\bar X_0, \mathbb{Q}_l)/\prod_{i=0}^n(1-q^it)$$ and $\det(1-F^*t, H^n(\bar X_0, \mathbb{Q}_l)$ is a polynomial with integer coefficients independent of $l$. 

Let's vary $X_0$ within the Lefschetz pencil of hypersurfaces that is defined over $\mathbb{F}_q$ (see (5.7) for $X=\mathbb{P}^{n+1}$; the existence of such a pencil is not clear; if we wanted to complete the argument sketched here, we would have to use the arguments that will be given in (7.1). One verifies that $E$ coincides here with the whole $H^n$ and (3.2) provides the Weil conjecture for all the hypersurfaces of the pencil, in particular for $X_0$.

\textbf{(5.13)} \textit{Bibliographical notes on paragraphs 3 and 4.}

A) The results of Lefschetz (4.1) and (5.1) to (5.5) are contained in this book [2]. For the local theory (4.1), it may be more handy to consult SGA 7, XIV (3.2).

B) The results of paragraph 4 are proved in parts XIII, XIV and XV of SGA 7.

C) (5.7) is proved in SGA 7, XVII\footnote{Also see the sketch of the proof presented in Milne.}.

D) (5.8) is proved in SGA 7, XVIII. The irreducibility theorem is proved there for $E$ but only under the hypothesis that $E \cap E^{\perp}=\{0\}$. Proof of the general case (for $E/(E \cap E^{\perp})$) is similar. 

\section{The rationality theorem}
\textbf{(6.1)} Let $\mathbb{P}_0$ be a projective space of dimension $\geq 1$ over $\mathbb{F}_q$, $X_0 \subset \mathbb{P}_0$ a projective nonsingular variety, $A_0 \subset \mathbb{P}_0$ a linear subspace of codimension two, $D_0 \subset  \mathbb{\widecheck P}_0$ the dual line, $\mathbb{\bar F}_q$ the algebraic closure of $\mathbb{F}_q$ and $\mathbb{P}, X, A, D$ over $\mathbb{\bar F}_q$ obtained from $\mathbb{P}_0, X_0, A_0, D_0$ by extension of scalars. The diagram (5.1.1) from (5.6) comes from a similar diagram over $\mathbb{F}_q$:
$$\begin{tikzcd}
X_0 \arrow[r, leftarrow, "\pi_0"]
& \tilde X_0 \arrow[d, "f_0"]\\
& D_0
\end{tikzcd} \ \ \textbf{(6.1.1)}$$

We assume that $X$ is connected of \textit{even} dimension $n+1=2m+2$ and that the pencil of hyperplane sections of $X$ defined by $D$ is a \textit{Lefschetz pencil}. The set $S$ of $t \in D$ such that $X_t$ is singular and defined over $\mathbb{F}_q$ comes from $S_0 \subset D_0$. We denote $U_0=D_0-S_0$ and $U=D-S$.

Let $u \in U$. The vanishing part of the cohomology $E \subset H^n(X_u, \mathbb{Q}_l)$ is stable under $\pi_1(U, u)$, so it is defined over $U$ by a local subsystem $\mathcal{E}$ of $R^nf_*\mathbb{Q}_l$. The latter is defined over $\mathbb{F}_q$: $R^if_*\mathbb{Q}_l$ is the inverse image of the $\mathbb{Q}_l$ sheaf $R^if_{0*}\mathbb{Q}_l$ on $D_0$ and, on $U$, $\mathcal{E}$ is the inverse image of a local subsystem $$\mathcal{E}_0 \subset R^nf_{0*}\mathbb{Q}_l.$$ The cup product is a skew-symmetric form $$\psi: R^nf_{0*}\mathbb{Q}_l \otimes R^nf_{0*}\mathbb{Q}_l \to \mathbb{Q}_l(-n).$$ Denoting by $\mathcal{E}_0^{\perp}$ the orthogonal of $\mathcal{E}_0$ relative to $\psi$, on $R^nf_{0*}\mathbb{Q}_l | U_0$ we see that $\psi$ induces a a perfect pairing $$\psi: \mathcal{E}_0/(\mathcal{E}_0 \cap \mathcal{E}_0^{\perp}) \otimes \mathcal{E}_0/(\mathcal{E}_0 \cap \mathcal{E}_0^{\perp}) \to \mathbb{Q}_l(-n).$$ 

\begin{thrm8}
For all $x \in |U_0|$ the polynomial $\det (1-F_x^*t, \mathcal{E}_0/(\mathcal{E}_0 \cap \mathcal{E}_0^{\perp}))$ has rational coefficients. 
\end{thrm8}

\begin{cor3}
Let $j_0$ be the inclusion of $U_0$ in $D_0$ and $j$ that of $U$ in $D$. The eigenvalues of $F^*$ acting on $H^1(D, j_*\mathcal{E}_0/(\mathcal{E}_0 \cap \mathcal{E}_0^{\perp}))$ are algebraic numbers all of which complex conjugates $\alpha$ satisfy $$q^{\frac{n+1}{2}-\frac{1}{2}} \leq |\alpha| \leq q^{\frac{n+1}{2}+\frac{1}{2}}.$$
\end{cor3}

According to (5.10) and (6.2), the hypotheses of (3.2) are in fact verified for $(U_0, \mathcal{E}_0/(\mathcal{E}_0 \cap \mathcal{E}_0^{\perp}), \psi)$ for $\beta=n$ and we apply (3.9). 

\begin{lemma6}
Let $\mathcal{G}_0$ be a locally constant $\mathbb{Q}_l$-sheaf on $U_0$ such that its inverse image $\mathcal{G}$ on $U$ is a constant sheaf. Then there exist units $\alpha_i$ in $\mathbb{\bar Q}_l$ such that for each $x \in |U_0|$ we have $$\det(1-F_x^*t, \mathcal{G}_0)=\prod_i (1-\alpha_i^{\deg(x)} t).$$
\end{lemma6}

The lemma expresses the fact that $\mathcal{G}_0$ is the inverse image of a sheaf on $Spec(\mathbb{F}_q)$, namely, its direct image on $Spec(\mathbb{F}_q)$ \footnote{Indeed, if $\varepsilon: U_o \to Spec \mathbb{F}_q$ is the canonical morphism, the assumption implies that $\varepsilon^*\varepsilon_* \mathcal{G}_0 \to \mathcal{G}_0$ becomes an isomorphism once we base change to $\mathbb{\bar F}_q$.}. The latter identifies with an $l$-adic representation $G_0$ of $Gal(\mathbb{\bar F}_q/\mathbb{F}_q)$ and we have\footnote{To establish the lemma, we also use the relationship between the stalks of the sheaf and those of its inverse/direct image (see Milne).} $$\det(1-Ft, G_0)=\prod_i (1-\alpha_i t).$$ 

Lemma (6.4) applies to $R^if_{0*} \mathbb{Q}_l \ (i \neq n)$, to $R^n f_{0*} \mathbb{Q}_l/\mathcal{E}_0$ and to $\mathcal{E}_0 \cap \mathcal{E}_0^{\perp}$\footnote{To see that the inverse images of the latter two are constant, use the equivalence of categories of locally constant sheaves and continuous representations of $\pi_1$ and the fact that the monodromy action acts by deforming the cohomology by the vanishing cycles.} . 

For $x \in |U_0|$ the fiber $X_x=f_0^{-1}(x)$ is a variety over the finite field $k(x)$. If $\bar x$ is a point of $U$ above $x$, $X_{\bar x}$ is obtained from $X_x$ by extension of scalars of $k(x)$ to the algebraic closure $k(\bar x)=\mathbb{\bar F}_q$ and $H^i(X_{\bar x}, \mathbb{Q}_l)$ is the stalk of $R^if_*\mathbb{Q}_l$ at $\bar x$. The formula (1.5.4) for the variety $X_x$ over $k(x)$ gives $$Z(X_x, t)=\prod_i \det(1-F_x^* t, R^i f_{0*} \mathbb{Q}_l)^{(-1)^{i+1}}$$ and $Z(X_x, t)$ is a product of $$Z^f=\det(1-F_x^* t, R^n f_{0*} \mathbb{Q}_l/\mathcal{E}_0)\det(1-F_x^* t, \mathcal{E}_0 \cap \mathcal{E}_0^{\perp}) \prod_{i \neq n} \det(1-F_x^* t, R^i f_{0*} \mathbb{Q}_l)^{(-1)^{i+1}}$$ and $$Z^m=\det(1-F_x^* t, \mathcal{E}_0/(\mathcal{E}_0 \cap \mathcal{E}_0^{\perp}).$$ Put $\mathcal{F}_0=\mathcal{E}_0/(\mathcal{E}_0 \cap \mathcal{E}_0^{\perp})$, $\mathcal{F}=\mathcal{E}/(\mathcal{E} \cap \mathcal{E}^{\perp})$ and apply (6.4) to the factors of $Z^f$. We find that there exist $l$-adic units $\alpha_i \ (1 \leq i 
\leq N)$ and $\beta_j \ (1 \leq j \leq M)$ in $\mathbb{\bar Q}_l$ such that for all $x \in |U_0|$ $$Z(X_x, t)=\frac{\underset{i} \prod (1-\alpha_i^{\deg(x)}t)}{\underset{j} \prod (1-\beta_j^{\deg(x)}t)} \det(1-F_x^*t, \mathcal{F}_0)$$ and in particular the right side is in $\mathbb{Q}(t)$. If some $\alpha_i$ coincides with a $\beta_j$, we can simultaneously delete $\alpha_i$ from the family of $\alpha$ and $\beta_j$ from the family of $\beta$. Therefore, we may and do assume that $\alpha_i \neq \beta_j$ for all $i$ and all $j$. 

\textbf{(6.5)} It suffices to prove that polynomials $\underset{i} \prod (1-\alpha_i t)$ and $\underset{j} \prod (1-\beta_j t)$ have rational coefficients, i.e. the family of $\alpha_i$ (resp. the family of $\beta_j$) is defined over $\mathbb{Q}$. We will deduce that from the following propositions. 

\begin{fact2}
Let $(\gamma_i) \ (1 \leq i \leq P)$ and $(\delta_j) \ (1 \leq j \leq Q)$ be two families of $l$-adic units in $\mathbb{\bar Q}_l$. Assume that $\gamma_i \neq \delta_j$. If $K$ is a large enough set of integers $\neq 1$, and $L$ is a large enough nowhere dense subset of $|U_0|$, then, if $x \in |U_0|$ satisfies $k \nmid \deg(x)$ (for all $k \in K$) and $x \notin L$, the denominator of $$\det(1-F_x^*t, \mathcal{F}_0) \underset{i} \prod (1-\gamma_i^{\deg(x)}t)/ \underset{j} \prod (1-\delta_j^{\deg(x)}t) \ \ \textbf{(6.6.1)}$$ written in irreducible form, is $\underset{j} \prod (1-\delta_j^{\deg(x)}t)$. 
\end{fact2}

The proof will be given in (6.10-13). According to (6.7) below, (6.6) provides an intrinsic description of the family of $\delta_j$ in terms of the family of rational fractions (6.6.1) for $x \in |U_0|$. 

\begin{lemma7}
Let $K$ be a finite set of integers $\neq 1$ and $(\delta_j) \ (1 \leq j \leq Q)$ and $(\varepsilon_j) \ (1 \leq j \leq Q)$ be two families of elements of a field. If, for all $n$ large enough, not divisible by any of the $k \in K$, the family of $\delta_j^n$ coincides with that of $\varepsilon_j^n$ (up to order), then the family of $\delta_j$ coincides with that of $\varepsilon_j$ (up to order). 
\end{lemma7}

We proceed by induction on $Q$. The set of integers $n$ such that $\delta_Q^n=\varepsilon_j^n$ is an ideal $(n_j)$. Let's prove that there exists a $j_0$ such that $\delta_Q=\varepsilon_{j_0}$. Otherwise the $n_j$ would be distinct from $1$ and there would be arbitrarily large integers $n$, not divisible by any of the $n_j$ nor by any of the $k \in K$. We would have $\delta_Q^n \neq \varepsilon_j^n$ and this contradicts the hypothesis. So there exists a $j_0$ such that $\delta_Q=\varepsilon_{j_0}$.We conclude by applying the induction hypothesis to the families $(\delta_j) \ (j \neq Q)$ and $(\varepsilon_j) \ (j \neq j_0)$.

\begin{fact3}
Let $(\gamma_i) \ (1 \leq i \leq P)$ and $(\delta_j) \ (1 \leq j \leq Q)$ be two families of $p$-adic units in $\mathbb{\bar Q}_l$, $R(t)=\underset{i} \prod (1-\gamma_i t)$ and $S(t)=\underset{j} \prod (1-\delta_j t)$. Assume that for all $x \in |U_0| \ \ \ \underset{j} \prod (1-\delta_j^{\deg(x)} t)$ \ \footnote{In the original paper $i$ appears instead of $j$.} divides $$\underset{i} \prod (1-\gamma_i^{\deg(x)}t) \det(1-F_x^*t, \mathcal{F}_0).$$ Then $S(t)$ divides $R(t)$. 
\end{fact3}
Remove from the families $(\gamma_i)$ and $(\delta_j)$ pairs of common elements until they verify the hypothesis of (6.6). Apply (6.6). By hypothesis, the rational fractions (6.6.1) are polynomials. Therefore, no $\delta$ survives, which means that\footnote{By the remarks after Proposition (6.6) and Lemma (6.7).} $S(t)$ divides $R(t)$. 

\textbf{(6.9)} We prove (6.5) and (6.2) (modulo (6.6)). Let's put $(\gamma_i)=(\alpha_i)$ and $(\delta_i)=(\beta_i)$ in (6.6). We get an intrinsic characterization of the family of $\beta_j$ in terms of the family of rational functions $Z(X_x, t) \ (x \in |U_0|)$. These being in $\mathbb{Q}(t)$, the family of $\beta_j$ is defined over $\mathbb{Q}$. 

Polynomials\footnote{for various $x$} $\underset{i} \prod (1-\alpha_j^{\deg(x)}t) \det (1-F_x^* t, \mathcal{F}_0)$ are therefore in $\mathbb{Q}[t]$. Proposition (6.8) provides an intrinsic description of the family of $\alpha_i$ in terms of this family of polynomials\footnote{Since a polynomial is uniquely characterized by its quotients.}. The family of $\alpha_i$ is thus defined over $\mathbb{Q}$. 

\textbf{(6.10)}
Let $u \in U$ and $\mathcal{F}_u$ the stalk of $\mathcal{F}$ at $u$. The arithmetic fundamental group $\pi_1(U_0, u)$, the extension of $\mathbb{\hat Z}=Gal( \mathbb{\bar F}/\mathbb{F})$ (generator: $\varphi$) by the geometric fundamental group\footnote{Both of the aforementioned groups are also called etale fundamental groups (of $U_0$ and $U$ respectively).} $\pi_1(U, u)$,  acts on $\mathcal{F}_u$ by symplectic similitudes\footnote{i.e. a group $A(V)$ of linear transformations $g$ of a vector space $V$ equipped with a nondegenerate bilinear form $<, >: V \times V \to k$ such that $<gv, gw>=\mu(g)<v, w>$ for a multiplicative character $\mu:A(V) \to k^{\times}$ called the similitude multiplier} $$\rho: \pi_1(U_0, u) \to CSp(\mathcal{F}_u, \psi).$$ We denote by $\mu(g)$ the multiplier of the symplectic similitude $g$. Let $$H \subset \mathbb{\hat Z} \times CSp(\mathcal{F}_u, \psi)$$ be the subgroup defined by the equation $$q^{-n}=\mu(g)$$ ($q$ being an $l$-adic unit, $q^n \in \mathbb{Q}_l^*$ is defined for all $n \in \mathbb{\hat Z}$). The fact that $\psi$ has values in $\mathbb{Q}_l(-n)$ implies that the map from $\pi_1$ to $\mathbb{\hat Z} \times CSp$, the coordinates of the canonical projection to $\mathbb{\hat Z}$ and $\rho$ factor through $$\rho_1: \pi_1(U_0, u) \to H.$$  
\begin{lemma8}
The image $H_1$ of $\rho_1$ is open in $H$. 
\end{lemma8}
Indeed, $\mathbb{\pi}_1(U_0, u)$ projects \textit{onto} $\mathbb{\hat Z}$ and the image of $\pi_1(U, u)=Ker(\pi_1(U_0, u) \to \mathbb{\hat Z})$ in \\ $Sp(\mathcal{F}_u, \psi)=Ker(H \to \mathbb{\hat Z})$ is open (5.10).   

\begin{lemma9}
For $\delta \in \mathbb{\bar Q}_l$ an $l$-adic unit, the set $Z$ of $(n, g) \in H_1$ such that $\delta^n$ is an eigenvalue of $g$ is closed of measure 0 \ \footnote{With respect to the Haar measure on $\mathbb{\hat Z} \times CSp$.}. 
\end{lemma9}
It is clear that $Z$ is closed. For each $n \in \mathbb{\hat Z}$ let $CSp_n$ be the set of $g \in CSp(\mathcal{F}_u, \psi)$ such that $\mu(g)=q^{-n}$ and let $Z_n$ be the set of $g \in CSp_n$ such that $\delta^n$ is an eigenvalue of $g$. Then $CSp_n$ is a homogeneous space for $Sp$ and we check that $Z_n$ is a proper algebraic subspace, thus, of measure $0$. According to (6.11), $H^1 \cap (\{n \} \times Z_n)$ is of measure $0$ in the inverse image in $H_1$ of $n$ and we apply Fubini to the projection $H_1 \to \mathbb{\hat Z}$. 

\textbf{(6.13)} Let us prove (6.6). For each $i$ and $j$, the set of integers $n$ such that $\gamma_i^n=\delta_j^n$ is the set of multiples of a fixed integer $n_{ij}$ (we do not exclude $n_{ij}=0$). By hypothesis, $n_{ij} \neq 1$. 

According to (6.12) and the Chebotarev's density theorem, the set of $x \in |U_0|$ such that $\beta_j^{\deg(x)}$ is an eigenvalue of $F_x^*$ acting on $\mathcal{F}_0$ is nowhere dense. We take for $K$ the set of $n_{ij}$ and for $L$ the set of $x$ as above. 
\section{Completion of the proof of (1.7)}

\begin{lemma10}
Let $X_0$ be a nonsingular absolutely irreducible\footnote{or geometrically irreducible} projective variety of even dimension over $\mathbb{F}_q$. Let $X$ over $\mathbb{\bar F}_q$ be obtained from $X_0$ by extension of scalars and $\alpha$ an eigenvalue of $F^*$ acting on $H^d(X, \mathbb{Q}_l)$. Then $\alpha$ is an algebraic number all of which complex conjugates, still denoted $\alpha$, satisfy $$q^{\frac{d}{2}-\frac{1}{2}} \leq |\alpha| \leq q^{\frac{d}{2}+\frac{1}{2}} \ \ \textbf{(7.1.1)}$$ 
\end{lemma10}
We proceed by induction on $d$ (always assumed even). The case $d=0$ is trivial even without assuming that $X_0$ is absolutely irreducible; we assume from now on $d \geq 2$. We put $d=n+1=2m+2$. 

If $\mathbb{F}_{q^r}$ is an extension of degree $r$ of $\mathbb{F}_q$ and $X_0'/\mathbb{F}_{q^r}$ is obtained from $X_0/\mathbb{F}_q$ by extension of scalars, the statement (7.1) for $X_0/\mathbb{F}_q$ is equivalent to (7.1) for $X_0'/\mathbb{F}_{q^r}$; in the same way as $q$ is replaced by $q^r$, the eigenvalues of $F^*$ are replaced by their $r$-th powers. 

According to (5.7), in a suitable projective embedding $i: X \to \mathbb{P}$, $X$ admits a Lefschetz pencil of hyperplane sections. The previous remark allows us to assume that the pencil is defined over $\mathbb{F}_q$ (once we replace $\mathbb{F}_q$ by a finite extension). 

Therefore, assume that there exists a projective embedding $X_0 \to \mathbb{P}_0$ and a subspace $A_0 \subset \mathbb{P}_0$ of codimension two that defines the Lefschetz pencil. We recall the notations of (6.1) and (6.3). A new extension of scalars allows us to assume that: 

a) The points of $S$ are defined over $\mathbb{F}_q$. 

b) The vanishing cycles for $x_s \ (s \in S)$ are defined over $\mathbb{F}_q$ (since only $\pm \delta$ is intrinsic, they can only be defined over quadratic extensions).

c) There exists a rational point $u_0 \in U_0$. We take the corresponding point $u$ of $U$ as the base point. 

d) $X_{u_0}=f_0^{-1}(u_0)$ admits a smooth hyperplane section $Y_0$ defined over $\mathbb{F}_q$. We let $Y=Y_0 \otimes_{\mathbb{F}_q} \mathbb{\bar F}_q$.

Since $\tilde X$ is obtained from $X$ by blowing up along a smooth subvariety $A \cap X$ of dimension two, we have $$H^i(X, \mathbb{Q}_l) \xhookrightarrow{} H^i(\tilde X, \mathbb{Q}_l)$$ (in fact, $H^i(\tilde X, \mathbb{Q}_l)= H^i(X, \mathbb{Q}_l) \oplus H^{i-2} (A \cap X, \mathbb{Q}_l)(-1))$\footnote{This is by the Thom isomorphism theorem. See a (rather technical) proof in Milne.}. It suffices to prove (7.1.1) for the eigenvalues $\alpha$ of $F^*$ acting on $H^d(\tilde X, \mathbb{Q}_l)$. 

The Leray spectral sequence for $f$ is $$E_2^{pq}=H^p(D, R^q f_*\mathbb{Q}_l) \Rightarrow H^{p+q}(\tilde X, \mathbb{Q}_l).$$ It suffices to prove (7.1.1) for the eigenvalues of $F^*$ acting on $E_2^{pq}$ for $p+q=d=n+1$. Those are\footnote{I highly recommend consulting Milne's book for the explanations of the steps in A), B), C).}:

A) $E_2^{2, n-1}$. According to (5.8), $R^{n-1}f_*\mathbb{Q}_l$ is constant. From (2.10) we have $$E_2^{2, n-1}=H^{n-1}(X_u, \mathbb{Q}_l)(-1).$$ Applying the weak Lefschetz theorem (corollary of SGA 4, XIV (3.2)) and the Poincare duality (SGA 4, XVIII), we have $$H^{n-1}(X_u, \mathbb{Q}_l)(-1) \xhookrightarrow{} H^{n-1}(Y, \mathbb{Q}_l)(-1)$$ and we apply the induction hypothesis to $Y_0$. 

B) $E_2^{0, n+1}$. If the vanishing cycles are nonzero, $R^{n+1}f_* \mathbb{Q}_l$ is constant and $$E_2^{0, n+1}=H^{n+1}(X_u, \mathbb{Q}_l).$$ The Gysin map $$H^{n-1}(Y, \mathbb{Q}_l)(-1) \to H^{n+1}(X_u, \mathbb{Q}_l)$$ is surjective (by an argument dual to that of A)) and we apply the induction hypothesis to $Y_0$. 

If the vanishing cycles are zero, the exact sequence of (5.8) b) gives the following exact sequence $$\underset{s \in S} \oplus \mathbb{Q}_l (m-n) \to E_2^{0, n+1} \to H^{n+1}(X_u, \mathbb{Q}_l).$$ The eigenvalues of $F$ acting on $\mathbb{Q}_l (m-n)$ are $q^{d/2}$ and for $H^{n+1}$ everything is as above. 

C) $E_2^{1, n}$. If we had the hard Lefschietz theorem, we would know that $\mathcal{E} \cap \mathcal{E}^{\perp}$ is zero and that $R^n f_* \mathbb{Q}_l$ is the direct sum of $j_*\mathcal{E}$ and a constant sheaf. The $H^1$ of a constant sheaf on $\mathbb{P}^1$ is zero and it would suffice to apply (6.3).

Since we have not proved the hard Lefshetz theorem yet, we will have to figure a way out. If the vanishing cycles are zero, $R^nf_* \mathbb{Q}_l$ is constant ((5.8) b)) and $E_2^{1, n}=0$. Therefore we may and do assume that the vanishing cycles are nonzero. Filter $R^nf_* \mathbb{Q}_l=j_*j^* R^n f_* \mathbb{Q}_l$ (5.8) by the subsheafs $j_*\mathcal{E}$ and $j_*(\mathcal{E} \cap \mathcal{E}^{\perp})$. If the vanishing cycles $\delta$ are not in $\mathcal{E} \cap \mathcal{E}^{\perp}$ \ \footnote{More precisely, not in $E \cap E^{\perp}$.} we have exact sequences\footnote{The constant sheaf in (7.1.2) is $j_*(R^n f_* \mathbb{Q}_l/ \mathcal{E})$. To see that this (and the next few) sheaves are constant, reason as in footnote 31 and use that the inverse/direct images of constant sheaves via $j$ are constant (the second claim is false in general).}: $$0 \to j_*\mathcal{E} \to R^n f_* \mathbb{Q}_l \to constant \ sheaf \to 0 \ \ \textbf{(7.1.2)}$$ $$0 \to \ constant \ sheaf \ j_*(\mathcal{E} \cap \mathcal{E}^{\perp})) \to j_*\mathcal{E} \to j_*(\mathcal{E}/(\mathcal{E} \cap \mathcal{E}^{\perp})) \to 0 \ \ \textbf{(7.1.3)}$$
If, God forbid, the $\delta$ are in $\mathcal{E} \cap \mathcal{E}^{\perp}$, we have $\mathcal{E} \subset \mathcal{E}^{\perp}$ and exact sequences\footnote{For $\mathcal{F}=R^nf_*\mathbb{Q}_l/j_* \mathcal{E}^{\perp}$.}: $$ 0 \to \ the \ constant \ sheaf \ j_*\mathcal{E}^{\perp} \to R^n f_* \mathbb{Q}_l \to \ a \ sheaf \ \mathcal{F} \to 0 \ \ \textbf{(7.1.4)}$$ $$ 0 \to \mathcal{F} \to \ the \ constant \ sheaf \ j_*j^* \mathcal{F} \to \underset{s \in S} \oplus \mathbb{Q}_l(n-m)_s \to 0 \ \ \textbf{(7.1.5)}$$ 

In the first case the long exact sequences in cohomology give $$H^1(D, j_* \mathcal{E}) \to H^1(D, R^n f_* \mathbb{Q}_l) \to 0 \ \ \textbf{(7.1.2')}$$ $$0 \to H^1(D, j_*\mathcal{E}) \to H^1(D, j_*(\mathcal{E}/(\mathcal{E} \cap \mathcal{E}^{\perp}))) \ \ \textbf{(7.1.3')}$$ and we apply (6.3).

In the second case, they give $$0 \to H^1(D, R^n f_*\mathbb{Q}_l) \to H^1(D, \mathcal{F}) \ \ \textbf{(7.1.4')}$$ $$ \underset{s \in S} \oplus \mathbb{Q}_l(n-m) \to H^1(D, \mathcal{F}) \to 0 \ \ \textbf{(7.1.5')}$$ and we remark that $F$ acts on $\mathbb{Q}_l(n-m)$ by multiplication by $q^{d/2}$.   

\begin{lemma11}
Let $X_0$ be a nonsingular projective absolutely irreducible variety of dimension $d$ over $\mathbb{F}_q$. Let $X$ over $\mathbb{\bar F}_q$ be obtained from $X_0$ by extension of scalars and $\alpha$ be an eigenvalue of $F^*$ acting on $H^d(X, \mathbb{Q}_l)$. Then $\alpha$ is an algebraic number all of which complex conjugates, still denoted $\alpha$,  satisfy $$|\alpha|=q^{\frac{d}{2}}.$$
\end{lemma11}

We first prove that (7.2) $\Rightarrow$ (1.7). For $X_0$ projective nonsingular over $\mathbb{\bar F}_q$ we have to prove the following statements: 

\textit{$W(X_{0}, i)$. Let $X$ be obtained from $X_0$ by extension of scalars of $\mathbb{F}_q$ to $\mathbb{\bar F}_q$. If $\alpha$ is an eigenvalue of $F^*$ acting on $H^i(X, \mathbb{Q}_l)$, then $\alpha$ is an algebraic number all of which complex conjugates, still denoted $\alpha$, satisfy $|\alpha|=q^{i/2}$.}

a) If $\mathbb{F}_{q^n}$ is an an extension of degree $n$ of $\mathbb{F}_q$ and $X_0'/\mathbb{F}_{q^n}$ is obtained from $X_0/\mathbb{F}_q$ by extension of scalars, then $W(X_0, i)$ is equivalent to $W(X_0', i)$: the extension of scalars replaces $\alpha$ by $\alpha^n$ and $q$ by $q^n$.

b) If $X_0$ is purely of dimension $n$, $W(X_0, i)$ is equivalent to $W(X_0, 2n-i)$; this follows from Poincare duality\footnote{If $\alpha$ is an eigenvalue of $F^*$ acting on $H^i(X, \mathbb{Q}_l)$, then $q^n/\alpha$ is an eigenvalue of $F^*$ acting on $H^{2n-i}(X, \mathbb{Q}_l)$. See Milne.}.

c) If $X_0$ is a an union of irreducible $X_0^{\alpha}$, $W(X_0, i)$ is equivalent to the collection of $W(X_0^{\alpha}, i)$.

d) If $X_0$ is purely of dimension $n$, $Y_0$ is a smooth hyperplane section of $X_0$ and $i < n$, then $W(Y_0, i) \Rightarrow W(X_0, i)$: this follows from the weak Lefschetz theorem\footnote{Depending on what one takes for the weak Lefschetz theorem, this follows either directly from it and the properties of Frobenius (if we assume that the weak Lefschetz theorem provides us with surjectivity of the Gysin map in cohomology) or one has to also invoke the Gysin sequence to establish surjectivity (if we only assume $H^i(X, \mathbb{Q}_l)=0$ for $i>d$ and $X$ affine). One can also avoid using the weak Lefschetz theorem and apply a Kunneth formula argument instead (due to A. Mellit). See Milne for clarifications.}. 

To prove the statements $W(X_0, i)$ we move in succession:

-by c) we assume that $X_0$ is purely of dimension $n$;

-by b) we also assume that $0 \leq i \leq n$;

-by a) and d) we also assume $i=n$;

-by a) and c) we also assume that $X_0$ is absolutely irreducible. 

Now the case satisfies the conditions of (7.2). 

\textbf{(7.3)} We prove (7.2). For every integer $k$, $\alpha^k$ is an eigenvalue of $F^*$ acting on $H^{kd}(X^k, \mathbb{Q}_l)$ (Kunneth's formula). For $k$ even, $X^k$ satisfies the conditions of (7.1), so we have $$q^{\frac{kd}{2}-\frac{1}{2}} \leq |\alpha^k| \leq q^{\frac{kd}{2}+\frac{1}{2}}$$ and $$q^{\frac{d}{2}-\frac{1}{2k}} \leq |\alpha| \leq q^{\frac{d}{2}+\frac{1}{2k}}.$$ Letting $k$ go to infinity, we establish (7.2). 
\section{First applications}
\begin{thrm9}Let $X_0 \subset \mathbb{P}_0^{n+r}$ be a nonsingular complete intersection over $\mathbb{F}_q$ of dimension $n$ and of multidegree $(d_1, \cdots d_r)$. Let $b'$ be the $n$-th Betti number of the complex nonsingular complete intersection with the same dimension and multidegree. Put $b=b'$ for $n$ odd and $b=b'-1$ for $n$ even. Then $$|\# X_0(\mathbb{F}_q)-\# \mathbb{P}^n(\mathbb{F}_q)| \leq b q^{n/2}.$$
\end{thrm9}

Let $X/\mathbb{\bar F}_q$ be obtained from $X_0$ and $\mathbb{Q}_l \eta^i$ be the line in $H^{2n}(X, \mathbb{Q}_l)$ generated by the $i$-th cup power of the cohomology class of the hyperplane section. On this line $F^*$ acts by multiplication by $q^i$. The cohomology of $X$ is the direct sum of the $\mathbb{Q}_l \eta^i \ (0 \leq i \leq n)$ and the primitive part of $H^n(X, \mathbb{Q}_l)$ of dimension $b$. According to (1.5), therefore, there exist $b$ algebraic numbers $\alpha_j$, the eigenvalues of $F^*$ acting on this primitive cohomology, such that $$\# X_0(\mathbb{F}_q)=\sum_{i=0}^n q^i+(-1)^n \sum_j \alpha_j.$$ According to (1.7), $|\alpha_j|=q^{n/2}$ and $$|\# X_0(\mathbb{F}_q)-\# \mathbb{P}^n(\mathbb{F}_q)|=|\# X_0(\mathbb{F}_q)-\sum_{i=0}^n q^i|= |\sum_j \alpha_j| \leq \sum_j |\alpha_j|=b q^{n/2}.$$

\begin{thrm10}
Let $N$ be an integer $\geq 1$, $\varepsilon: (\mathbb{Z}/N)^* \to \mathbb{C}^*$ a character, $k$ an integer $\geq 2$ and $f$ a holomorphic modular form on $\Gamma_0(N)$ of weight $k$ and with character $\varepsilon$ : $f$ is  a holomorphic function on the the Poincare half-plane $X$ such that for $\left(\begin{array}{cc} a & b\\ c & d \end{array}\right) \in SL(2, \mathbb{Z}), $ with $c \equiv 0 \ (N)$ we have $$f \left(\frac{az+b}{cz+d} \right)=\varepsilon(a)^{-1}(cz+d)^kf(z).$$ We assume that $f$ is cuspidal and primitive ("new" in the sence of Arkin-Lehner and Miyake), in particular an eigenvector of the Hecke operators $T_p \ (p \nmid N)$. Let $f=\sum_{n=1}^{\infty} a_nq^n$ with $q=e^{2\pi i z}$ (and $a_1=1$). Then for $p$ prime not dividing $N$ $$|a_p| \leq 2p^{\frac{k-1}{2}}.$$ In other words, the roots of the equation $$T^2-a_p T+\varepsilon(p)p^{k-1}$$ are of absolute value $p^{\frac{k-1}{2}}$. 
\end{thrm10}

These roots are indeed the eigenvalues of the Frobenius acting on $H^{k-1}$ of a nonsingular projective variety of dimension $k-1$ defined over $\mathbb{F}_p$.

Under restrictive assumptions, this fact is proved in my Bourbaki expose (Formes modulaires et representations $l$-adiques, expose 355, February 1969, in: \textit{Lecture Notes in Mathematics, 179)}. The general case is not much more difficult. 

\textit{Remark} \textbf{(8.3)} J.P.Serre and myself have recently proved that (8.2) remains true for $k=1$. The proof is quite different. 

The following application was suggested to me by E.Bombiery.
\begin{thrm11}
Let $Q$ be a polynomial in $n$ variables and of degree $d$ over $\mathbb{F}_q$, $Q_d$ a homogeneous part of degree $d$ of $Q$ and $\psi: \mathbb{F}_q \to \mathbb{C}^*$ an additive nontrivial character on $\mathbb{F}_q$. We assume that: 

(i) $d$ is coprime to $p$ \footnote{(the characteristic of $\mathbb{F}_q)$}

(ii) The hypersurface $H_0$ in $\mathbb{P}_{\mathbb{F}_q}^{n-1}$ defined by $Q_d$ is smooth. 

Then $$|\sum_{x_1, \cdots, x_n \in \mathbb{F}_q} \psi(Q(x_1, \cdots, x_n))| \leq (d-1)^n q^{n/2}.$$
\end{thrm11}

After replacing $Q$ by a scalar multiple, we may (and do) assume that\footnote{For $q=p^f$ we have $Tr_{\mathbb{F}_q/\mathbb{F}_p}(x)=x+x^p+\cdots+x^{p^{f-1}}$.} $$\psi(x)=\exp(2 \pi i Tr_{\mathbb{F}_q/\mathbb{F}_p}(x)/p) \ \ \textbf{(8.4.1)}$$ Let $X_0$ be an etale covering of the affine space $\mathbb{A}_0$ of dimension $n$ over $\mathbb{F}_q$ with equation $T^p-T=Q$ and let $\sigma$ be the projection of $X_0$ to $\mathbb{A}_0$: $$\sigma: X_0 \to \mathbb{A}_0$$ $$X_0=Spec(\mathbb{F}_q [x_1, \cdots, x_n, T]/(T^p-T-Q)).$$ The covering $X_0$ is Galois with Galois group $\mathbb{Z}/p$; $i \in \mathbb{Z}/p=\mathbb{F}_p$ acts by $T \to T+i$. 

We let $x \in \mathbb{A}_0(\mathbb{F}_q)$ and compute the Frobenius endomorphism on the fiber of $X_0/\mathbb{A}_0$ at $x$. Let $q=p^f$ and let $\mathbb{\bar F}_q$ be the algebraic closure of $\mathbb{F}_q$. For $(x, T) \in X_0(\mathbb{\bar F})$ above $x$ we have $F((x, T))=(x, T^q)$ and $$T^q=T+\sum_{i=1}^f(T^{p^i}-T^{p^{i-1}})=T+\sum Q(x)^{p^{i-1}}=T+Tr_{\mathbb{F}_q/\mathbb{F}_p}(Q(x)).$$ This is the action of the element $Tr_{\mathbb{F}_q/\mathbb{F}_p}(Q(x))$ of the Galois group. 

Let $E$ be the field of the $p$-th roots of unity and $\lambda$ a finite place of $E$ coprime to $p$. We will work in $\lambda$-adic cohomology. For $j \in \mathbb{Z}/p$, let $\mathcal{F}_{j, 0}$ be a $E_{\lambda}$ local system of rank one on $\mathbb{A}_0$ defined by $X_0$ and $\psi(-jx): \mathbb{Z}/p \to E^* \to E_{\lambda}^*$: we have $\iota:X_0 \to \mathcal{F}_{j, 0}$ and $\iota(i \star x)=\psi(-ij)\iota(x)$. 
Denote without $_0$ objects obtained from $\mathbb{A}_0, X_0, \mathcal{F}_{j, 0}$ by extension of scalars to $\mathbb{\bar F}_q$. The trace formula (1.12.1) for $\mathcal{F}_{j, 0}$ gives: $$\sum_{x_1, \cdots, x_n \in \mathbb{F}_q} \psi(Q(x_1, \cdots, x_n))=\sum_i Tr(F^*, H_c^i(\mathbb{A}, \mathcal{F}_1)) \ \ \textbf{(8.4.2)}$$ 

We have $\sigma_* E_{\lambda}=\underset{j} \oplus \mathcal{F}_j$ and so $$H_c^*(X, \mathbb{Q}_l) \otimes_{\mathbb{Q}_l} E_{\lambda}=\underset{j} \oplus H_c^*(\mathbb{A}, \mathcal{F}_j) \ \ \textbf{(8.4.3)}$$ For $j=0$, $\mathcal{F}_j$ is the constant sheaf $E_\lambda$; this factor corresponds to inclusion, by taking the inverse image, of the cohomology of $\mathbb{A}$ in that of $X$.

\begin{lemma12}
(i) For $j \neq 0$, $H_c^i(\mathbb{A}, \mathcal{F}_j)$ is zero for $i \neq n$, for $i=n$, the cohomology space has dimension $(d-1)^n$. 

(ii) For $j \neq 0$, the cup product $$H_c^n(\mathbb{A}, \mathcal{F}_j) \otimes H_c^n(\mathbb{A}, \mathcal{F}_{-j}) \to H^{2n}(\mathbb{A}, E_{\lambda}) \overset{Tr} \to E_{\lambda}(-n)$$ is a perfect pairing. 

(iii) $X_0$ is open in a nonsinqular projective variety $Z_0$. 
\end{lemma12}

Let's deduce (8.4) from (8.5). Let $j_0: X_0 \xhookrightarrow{} Z_0$ and $j: X \xhookrightarrow{} Z$ be obtained by extension of scalars of $\mathbb{\bar F}_q$. According to (8.4.2), (i) and (1.7) for $Z_0$, it suffices to prove the injectivity of $$H_c^n(\mathbb{A}, \mathcal{F}_1) \overset{\sigma^*} \to H_c^n(X, \mathcal{F}_1)=H_c^n(X, E_\lambda) \overset{j_!} \to H^n(Z, E_{\lambda}).$$ We have $Tr(a \cup b)=\frac{1}{p}Tr(j_! \sigma^* a \cap j_! \sigma^* b)$, so injectivity follows from (ii)\footnote{By the same reasoning as in the proof of (7.1) A) (for which I referred to Milne).}. 

\textbf{(8.6)} We prove (8.5) (iii). Let $\mathbb{P}_0$ be the projective space over $\mathbb{F}_q$, obtained from $\mathbb{A}_0$ by adding a hyperplane at infinity $\mathbb{P}_0^{\infty}$, $H_0 \subset \mathbb{P}_0^{\infty}$ with equation $Q_d=0$ and $Y_0$ the covering of $\mathbb{P}_0$ normalizing $\mathbb{P}_0$ along $X_0$. 
$$\begin{tikzcd}
X_0 \arrow[r, hook] \arrow[d, "\sigma"]
&  Y_0 \arrow[d]\\
\mathbb{A}_0 \arrow[r, hook] & \mathbb{P}_0 \arrow[r, hookleftarrow] & \mathbb{P}_0^{\infty} \arrow[r, hookleftarrow] & H_0
\end{tikzcd} \ \ \textbf{(8.7.1)}$$
Let's study $Y_0/\mathbb{P}_0$ near the infinity, locally for the etale topology. 

\begin{lemma13}
$Y_0$ is smooth outside the inverse image of $H_0$.
\end{lemma13}

The divisor of a rational function $Q$ on $\mathbb{P}_0$ is the sum of the finite part $div(Q)_f$ and $(-d)$ times the hyperplane at infinity. We have: $$div(Q)=div(Q)_f-d\mathbb{P}_0^{\infty} \ \ \textbf{(8.7.1)}$$ $$div(Q)_f \cap \mathbb{P}_0^{\infty}=H_0$$ At a finite distance, $Y_0=X_0$ is etale over $\mathbb{A}_0$, so smooth. At the infinity but outside the inverse image of $H_0$ there exist local coordinates $(z_1, \cdots, z_n)$ such that $Q=z_1^{-d}$ (here we use $(d, p)=1$). In these coordinates, $Y_0$ appears as a product of a curve and a smooth space (corresponding to coordinates $z_2, \cdots, z_n$). By normality it is smooth.  

\begin{lemma14}
In the etale neighborhood of a point above $H_0$, $Y_0$ is smooth on a normal singular surface, always the same. 
\end{lemma14}

This time we can find local coordinates such that $Q=z_1^{-d}z_2$. Indeed, since $H_0$ is smooth, $div(Q)_f$ is smooth in the neighborhood of infinity and crosses $\mathbb{P}_0^{\infty}$ transversely. This form is independent of the chosen point and uses only two coordinates, hence the assertion.

\textbf{(8.9)} The following method (due to Zariski) allows one to resolve singularities on surfaces: alternately, we normalize and we blow up the (reduced) singular locus. Operators in play commute with etale localization and taking products by a smooth space. The method of Zariski, therefore, allows one to resolve singularities on a space that (like $Y_0$) is, locally for the etale topology, smooth on a surface. The resolution obtained from $Y_0$ is the $Z_0$ we seek. 

If $T$ is a curve on a surface $S$ containing the singular locus and $T'$ is the inverse image of $T$ in the Zariski resolution $S'$ of $S$, we know that if we repeatedly blow up the (reduced) singular locus of $(T')_{red}$ in $S'$, we obtain a surface $S''$ such that the inverse reduced image $(T'')_{red}$ of $T$ in $S''$ is a divisor with normal crossings. Again, operations in play commute with etale localization and taking products by a smooth space. Reasoning as above and observing that $(Y_0, infinity)$ is locally smooth in $(S, T)$, we can find $Z_0$ such that $Z_0-X_0$ is a divisor with normal crossings.

\textbf{(8.10)} We prove (8.5) (i), (ii). These assertions are geometric; this allows us to work from now on in $\mathbb{\bar F}_q$. Let $S'$ be the affine space over $\mathbb{\bar F}_q$ that parametrizes polynomials in $n$ variables of degree $\leq d$ and let $S$ be an open in $S'$ corresponding to the polynomials, which homogeneous part of degree $d$ has nonzero discriminant. We denote by $Q_{S} \in H^0(S, \mathcal{O}[x_1, \cdots, x_n])$ the universal polynomial\footnote{$\sum_{i_1+ \cdots+ i_n=d}a_{i_1, \cdots, i_n}x_1^{i_1} \cdots x_n^{i_n}$} of $S$ and by $X_S$ the Galois etale covering of $\mathbb{A}_S=\mathbb{A}^n \times S$ with equation $T^p-T=Q_S$ and Galois group $\mathbb{Z}/p$. Let $\mathbb{P}_S=\mathbb{P}^n \times S$ be the projective completion of $\mathbb{A}_S$ and $Y_S$ the normalization of $\mathbb{P}_S$ along $X_s$. We have, for $S$, a diagram similar to (8.6.1). 

Expressions of $Q$ in local coordinates (8.7) and (8.8) remain valid in this situation with parameters such that locally for the etale topology on $Y_S$, $Y_S/S$ is isomorphic to the product of $S$ (that is smooth) with one fiber. The method of canonical resolution used in (8.9) gives us a relative compactification $Z_S/S$ of $X_S/S$ with $Z_S-X_S$ a divisor with normal crossings relative to $S$ 
$$\begin{tikzcd}
X_S \arrow[r, hook, "u"] \arrow[d, "\sigma"]
&  Z_S \arrow[d, "f"]\\
\mathbb{A}_S \arrow[r, hook, "a"] & S 
\end{tikzcd}$$ ($f$ is proper and smooth, $u$ is an open immersion, $Z_s-X_s$ a divisor with relative normal crossings).

Let $\mathcal{F}_{j, S}$ be a $E_{\lambda}$-sheaf on $\mathbb{A}_S$ obtained as in (8.4) from $X_S/\mathbb{A}_S$. We have $\sigma^*E_{\lambda}=\oplus \mathcal{F}_{j, S}$, so $$R^*(fu)_! (E_{\lambda})=\underset{j} \oplus R^* a_! \mathcal{F}_{j, S}.$$ The properties of $Z_S$ ensure that $R^i(fu)_!E_{\lambda}=R^if_*(u_!E_{\lambda})$ is a locally constant sheaf on $S$. 
Therefore, $R^ia_!\mathcal{F}_{j, S}$ is also locally constant.  Since $S$ is connected, it suffices to prove (8.5) (i), (ii) for a particular polynomial $Q$. We will take $Q=\sum_i x_i^d$. This polynomial satisfies the nonsingularity condition because $(d, p)=1$. For this polynomial variables in the exponential sum (8.4) are separated. This corresponds to the fact that $\mathcal{F}_j$ is the tensor product of the inverse images of similar sheaves $\mathcal{F}_j^1$ on the factors of dimension one $\mathbb{A}^1$ of $\mathbb{A}=\mathbb{A}^n$. By Kunneth's formula $$H^*(\mathbb{A}, \mathcal{F}_j)=\oplus H^*(\mathbb{A}^1, \mathcal{F}_j^1).$$ This reduces the proof of (8.5) (i), (ii) to the case when $n=1$ and $Q$ is $x^d$.

\textbf{(8.11)} Let's deal with this particular case. The covering $X$ of $\mathbb{A}$ is irreducible, so for $i=0, 2$ $$H_c^i(\mathbb{A}, E_{\lambda}) \overset{\sim} \to H_c^i(X, E_{\lambda}).$$ So for $i \neq 1$ and $j \neq 0$ we have\footnote{In the original paper $\mathbb{A}$ appears instead of $X$.} $$H_c^i(X, \mathcal{F}_j)=0.$$ Assertion (ii) follows from (2.8) or (2.12) and the fact that $u_!\mathcal{F}_j=u_*\mathcal{F}_j$. To prove (i) it remains to verify that $$\chi_c(\mathbb{A}, \mathcal{F}_j)=1-d.$$ According to the Euler-Poincare formula (see expose Bourbaki 286 of February 1965, by M.Raynaud), it is equivalent to the following lemma.

\begin{lemma15}
Swan's conductor of $\mathcal{F}_j$ at infinity equals $d$.
\end{lemma15}
This statement is equivalent to the following. 

\begin{lemma16}
Let $k$ be a finite field of characteristic $p$, $y \in k[[x]]$ an element of valuation $d$ coprime to $p$, $L$ the extension of $K=k((x))$ generated by the roots of $T^p-T=y^{-1}$ and $\chi$ the following character on $Gal(L/K)$ with values in $\mathbb{Z}/p$: $$\chi(\sigma)=\sigma T-T.$$ Then $\chi$ has conductor $d+1$. 
\end{lemma16}

By extension of the residue field we may assume that $k$ is algebraically closed rather than finite and apply: J.P.Serre, Sur les corps locaux a corps residuel algebriquement clos, \textit{Bull. Soc. Math. France}, 89 (1961), p. 105-154, $n^o$ 4.4.

\begin{flushright}
\textit{Manuscript received on 20 September 1973.}
\end{flushright}

\end{document}